\documentclass[12pt,
               reqno,
               oneside]{amsart} %

\usepackage[left=1.25in,right=1.25in,bottom=.7in]{geometry} 

\makeatletter
\renewenvironment{proof}[1][\proofname]{\par
  \pushQED{\qed}%
  \normalfont 
  \trivlist
  \item[\hskip\labelsep
        \itshape
    #1\@addpunct{.}]\ignorespaces
}{%
  \popQED\endtrivlist\@endpefalse
}
\makeatother

\usepackage{%
            amsfonts,
            amssymb,
            amsmath,
            amsthm,
            mdwlist%
}
\usepackage[mathscr]{eucal}

\renewcommand\qedsymbol{\ensuremath{\overstrike\sqcap\sqcup}}

\usepackage{graphicx}

\usepackage[oldenum]{paralist} 

\newcommand\thmlistfix{\vspace{-1.0\topsep}}

\usepackage{hyperref} 
\usepackage[all]{hypcap} 

\newtheorem{lem}{Lemma}            
\newtheorem{thm}{Theorem}          
\newtheorem{prop}[thm]{Proposition}
\newtheorem{cor}{Corollary}[thm]   
 
\newtheorem*{thm*}{Theorem}
\newtheorem{prop*}{Proposition}
\newtheorem*{cor*}{Corollary}

\newtheorem*{wrapping genera estimate}{Wrapping Genera Estimate}%
\newtheorem*{braid index estimate}{Braid Index Estimate}%
\newtheorem*{crossing number estimate}{Crossing Number Estimate}%
\newtheorem*{free genus estimate}{Free Genus Estimate}%

\theoremstyle{definition}
\newtheorem*{definition}{Definition}
\newtheorem*{definitions}{Definitions}
\newtheorem*{remark}{Remark}
\newtheorem*{remarks}{Remarks}
\newtheorem*{question}{Question}
\newtheorem*{questions}{Questions}
\newtheorem*{hist}{Historical note}
\newtheorem*{hists}{Historical notes}
\newtheorem*{notation*}{Notation}
\newtheorem*{convention}{Convention}

\newtheoremstyle{construction}%
  {3pt}
  {3pt}
  {}
  {}
  {\bfseries\upshape}
  {.}
  {.5em}
  {\thmname{#1}\thmnumber{ #2}\thmnote{ (#3)}}

\theoremstyle{construction}
\newtheorem{construction}{Construction}

\newenvironment{unproved}[1]{\def\QQQQ{\string#1}\begin{\QQQQ}}%
                                  {\hfill\qedsymbol\end{\QQQQ}}

\newcounter{xampl}
\newtheorem{example}[xampl]{Example}

\renewcommand\a{\alpha}
\renewcommand\b{\beta} 
\newcommand\e{\varepsilon} 
\newcommand\g{\gamma} 
\renewcommand\d{\delta} 
 
\renewcommand\phi{\varphi} 
\newcommand\s{\sigma} 
\renewcommand\t{\tau}
\newcommand\x{\xi}
\newcommand\D{\Delta}

\newcommand\overstrike[2]{{\setbox0\hbox{$#2$}\hbox to \wd0{\hss
                         $#1$\hss}\kern-\wd0\box0}}
\newcommand\twobars[5]{\vcenter{\hrule height#1ex width#2ex
                               \vskip#3ex
                               \hrule height#4ex width#5ex}}
\newcommand\stroke[3]{\vrule height#1ex width#2ex depth#3ex}
\newcommand\connsum{\mathop{\raise.25ex\hbox{\overstrike\parallel=}}}
\newcommand\Connsum{\mathop{\displaystyle\raise.25ex\hbox{\overstrike\parallel=}}}
\newcommand\bdconnsum{\mathop{{}\hskip.3ex
                   \stroke{1.7}{.06}{-.13}
                   \twobars{.06}{.7}{.8}{.06}{.7}
                   \stroke{1.03}{.06}{.54}
                   \hskip.3ex{}}}

\newcommand\plumb[1]{\ast_{#1}}

\makeatletter
\DeclareRobustCommand{\prescript}[3]{%
  \@mathmeasure\z@\scriptstyle{#1}
  \@mathmeasure\tw@\scriptstyle{#2}
  \ifdim\wd\tw@>\wd\z@
    \setbox\z@\hbox to\wd\tw@{\hfil\unhbox\z@}%
  \else
    \setbox\tw@\hbox to\wd\z@{\hfil\unhbox\tw@}%
  \fi
  \mathop{}%
  \mathopen{\vphantom{#3}}^{\box\z@}_{\box\tw@}%
  #3%
}
\makeatother

\newcommand{\after}{\mathop{\circ}} 
\newcommand\AKn[2]{A(#1,#2)}  
\renewcommand\arg{\operatorname{arg}}

\newcommand\Bd{\partial} 
\newcommand\bi{\operatorname{brin}} 
\newcommand\blowup{\Xi} 
\newcommand\bydeffont{\textit} 
\newcommand\bydeffontname{italic} 
\newcommand\bydef[1]{\bydeffont{#1}}

\newcommand\card[1]{\operatorname{card}(#1)} 
\newcommand\C{\mathbb C} 
\newcommand\Cext{\P_1(\C)} 
\newcommand\clint[2]{[#1,#2]} 
\newcommand\clopint[2]{[#1,#2{[}} 
\newcommand\conjug[1]{\overline{#1}} 
\newcommand\crit[1]{\operatorname{crit}(#1)} 
\newcommand\cross[1]{#1^{\scriptscriptstyle{\times}}} 

\newcommand\cut[2]{#1\,\wr\,#2} 

\newcommand\diff[2]{\mbox{${#1}\setminus{#2}$}} 
\newcommand\divisor[1]{\mathbb{D}(#1)} 
\newcommand\double[3]{D(#1,#2,#3)}

\renewcommand\emptyset{\varnothing} 
\newcommand\Ext[2]{\operatorname{Ext}({#2}\hookrightarrow{#1})} 
\newcommand\ext{\operatorname{ext}} 

\newcommand\foreignfont{\textit} 
\newcommand\foreign[1]{\foreignfont{#1}}
\newcommand\framing[1]{\mathbf{#1}}
\newcommand\freegenus[1]{g_{\textup{f}}(#1)} 
\newcommand\freerank[1]{b_{\textup{f}}(#1)} %
\newcommand\from{\colon\thinspace} 

\newcommand\G{\mathbb{G}}

\newcommand\Hs{\mathscr{H}}
\newcommand\h[2]{h^{\scriptscriptstyle(#1)}_{#2}} 
\newcommand\hb[1]{\Hs^{#1}} 

\newcommand\id[1]{\operatorname{id}_{#1}} 
\renewcommand\Im{\operatorname{Im}} 
\newcommand\imunit{i}
\newcommand\ind[2]{\operatorname{ind}(#1;#2)} 
\newcommand\Int{\operatorname{Int}} 
\newcommand\isdefinedas{\mathrel{:=}}
\newcommand\whichdefines{\mathrel{=:}}
\newcommand\iso{\linebreak[1]\cong\linebreak[1]} 

\newcommand\liftof[1]{\phantom{.}{\tilde{\vphantom{o}}\kern0.1ex}{#1}}

\newcommand\meromapto{\dashrightarrow}
\renewcommand\meromapto{\rightsquigarrow}

\newcommand\Mgood{\mbox{M-M--good}} 
\newcommand\Mmap[1]{\phi_{#1}} 
\newcommand\Mmapr[2]{\phi_{#1,#2}} 
\newcommand\MN{\mathscr{MN}} 
\newcommand\MNfree{\mathscr{MN}_\textup{f}} 

\newcommand\Nb[2]{\operatorname{Nb}({#2}\hookrightarrow{#1})} 
\newcommand\N{\mathbb N} 

\newcommand\opint[2]{{]}#1,#2{[}} 

\renewcommand\P{\mathbb P} 
\newcommand\page[2]{S(#1,#2)}
\newcommand\pr{\operatorname{pr}}

\newcommand\R{\mathbb R} 
\newcommand\rank{\operatorname{rank}}
\renewcommand\Re{\operatorname{Re}} 
\newcommand\reduced[1]{\widetilde#1} 
\newcommand\restr{{\smash[bt]{\mathrel\big|}}} 

\newcommand\Seif[1]{\begin{bmatrix}#1\end{bmatrix}} 
\newcommand\splicedlink[3]{\Psi(#1,#2,#3)}
\newcommand\splicedmap[3]{\Phi(#1,#2,#3)} 
\newcommand\sub{\subset} 
\newcommand\Suchthat{\thinspace\colon\thinspace} 
\newcommand\suchthat{\thinspace\colon\thinspace} 

\newcommand\multitwist[3]{\operatorname{\varsigma}_{#3}(#1,#2)} %

\newcommand\wrapgenus[1]{g_{\textup{wr}}(#1)}
\newcommand\wlapgenus[1]{g_{\textup{w}\ell}(#1)}

\newcommand\Z{\mathbb Z} 

\begin{document}

\bibliographystyle{amsplain}

\title[Bounds on the Morse--Novikov number] 
{Constructions of Morse maps for knots and links,\\
and upper bounds on the Morse--Novikov number}

\author{Mikami Hirasawa}
\address{Department of Mathematics, 
Gakushuin University, 1-5-1 Mejiro, Tosima-ku, Tokyo 171-8588 Japan}  
\email{\href{mailto:hirasawa@math.gakushuin.ac.jp}%
{hirasawa@math.gakushuin.ac.jp}}

\author{Lee Rudolph}
\address{Department of Mathematics and Computer Science, 
Clark University, 
Worcester MA 01610 USA}  
\email{\href{mailto:lrudolph@black.clarku.edu}%
{lrudolph@black.clarku.edu}}
\thanks{Partial support for both authors was supplied by
the Fonds National Suisse,
and for the first author by 
Fellowships of the Japan Society for the Promotion of 
Science for Japanese Junior Scientists
and MEXT Grant-in-Aid for Young Scientists (B) 14740048.}

\begin{abstract} 
The {Morse--Novikov number} $\MN(L)$ of an oriented link
$L \sub S^3$ is the minimum number of critical points 
of a Morse map $\diff{S^3}{L} \to S^1$ 
representing the class 
of a Seifert surface for $L$ in 
$H_1(S^3, L; \Z)$ 
(e.g., $\MN(L)=0$ if and only if $L$ is fibered). 
We develop various constructions of Morse maps
(Milnor maps, 
Stallings twists,
splicing along a link which is a closed braid
with respect to a Morse map,
Murasugi sums,
cutting a Morse map along an arc on a page)
and use them to bound Morse--Novikov numbers from above 
in terms of other knot and link
invariants (free genus, 
crossing number,
braid index,
wrapping genus
and layered wrapping genus).
\end{abstract}

\maketitle

\section{Introduction; statement of results}\label{introduction}

An oriented link $L\sub S^3$ determines 
a cohomology class 
$\xi_L\in H^1(\diff{S^3}{L};\Z)
\iso 
\pi_0(\operatorname{Map}(\diff{S^3}{L},K(\Z,1)))$.
The homotopy class of maps $\diff{S^3}{L}\to S^1=K(\Z,1)$ 
corresponding to $\xi_L$ contains smooth maps which are 
{Morse} (that is, have no degenerate critical points),
and which restrict to a standard fibration in a 
neighborhood of $L$ (so they have only finitely many 
critical points).  The minimum number of critical points 
of such a map is the Morse--Novikov number $\MN(L)$.  
Tautologously, $\MN(L)=0$ 
if and only if $L$ is a fibered link.  It is natural to ask 
how to calculate, or estimate, $\MN(L)$ for general $L$.
Moreover, for fibered links, there exist both 
nice characterizations in other terms (e.g., 
$\MN(L)=0$ if and only if the kernel of 
$\pi_1(\diff{S^3}{L})\to \Z\suchthat [\g]\mapsto \int_\g \xi_L$ 
is finitely generated)
and an array of interesting constructions (e.g.,
links of singularities and Milnor fibrations \cite{Milnor:singular-points}, 
Murasugi sums \cite{Murasugi:plumbing,Stallings,Gabai:Murasugi1},
Stallings twists \cite{Stallings},
splicing \cite{Eisenbud-Neumann}).
Again, it is natural to ask what happens 
in general.

Some progress on these questions was made in \cite{P-R-W}. 
There, the Morse--Novikov theory of maps from manifolds 
to the circle (introduced by Novikov \cite{Novikov:multivalued}, 
and previously applied to knots in $S^3$ by Lazarev \cite{Lazarev})
was applied to give lower bounds for $\MN(L)$: for example,
it was shown (using analogues for Novikov homology of the 
Morse inequalities for ordinary homology) that for all $n\ge 0$, 
there exists a knot $K$ with genus $g(K)=n$ and $\MN(K)\ge 2n$; 
it was also shown that there are knots with vanishing Novikov 
homology and non-zero Morse--Novikov number.
Subadditivity 
\begin{equation*}
\MN(L_0\connsum L_1)\leq \MN(L_0)+\MN(L_1).
\tag{$\sharp$}\label{subadditivity over connected sum:intro}
\end{equation*}
of Morse--Novikov number over connected sum
was established by an explicit construction, and
it was conjectured that a restatement 
\thetag{\ref{subadditivity over Murasugi sum:proof}}
of \thetag{\ref{subadditivity over connected sum:intro}}
in terms of Seifert surfaces extends to arbitrary Murasugi sums.  
(In fact, earlier work of Goda \cite{Goda:suture-sum,Goda:handle-number}
immediately implies \thetag{\ref{subadditivity over Murasugi sum:proof}};
Goda's results are stated in terms of his ``handle number'' 
of a Seifert surface $R$, not the Morse--Novikov number 
of the link $\Bd R$, and his exposition and proofs 
use Gabai's language of sutured manifolds \cite{Gabai:foliations-genera} 
and $C$-product decompositions \cite{Gabai:detect}, 
not that of Morse maps.)

The present paper continues the investigations of \cite{P-R-W}.
Section~\ref{preliminaries} assembles preliminary material
and generalities on Morse maps.
Section~\ref{Morse maps from Milnor maps} constructs
\hyperlink{Milnor maps hyperconstruction}{Milnor maps}
and provides two simple, but fundamental, 
examples of Milnor maps which are Morse maps but not
fibrations: $u$, a minimal Morse map for the (non-fibered)
$2$-component unlink $U$, and $o_1$, a Morse map with
two critical points for the (fibered) unknot $O$.
Section~\ref{twisting} describes 
\hyperlink{Stallings twists hyperconstruction}{Stallings twists} of 
Morse maps.
Section~\ref{monodromies} shows how 
to a Morse map $f\from \diff{S^3}{L}\to S^1$ 
(satisfying a condition much weaker than being a fibration) 
are associated certain maps of 
surfaces---\hyperlink{monodromy hyperconstruction}{monodromies
and adiexodons} (the latter being trivial when, and only when,
$f$ is a fibration)---from which $f$ can be reconstructed.
Section~\ref{splicing} introduces closed $f$-braids 
in $\diff{S^3}{L}$, and constructs Morse maps by 
\hyperlink{splicing hyperconstruction}{splicing} along 
closed $f$-braids.
Section~\ref{Murasugi sums} constructs Morse maps as 
\hyperlink{M. sum of Morse maps hyperconstruction}{Murasugi sums} 
of simpler Morse maps, and in particular 
by \hyperlink{cutting hyperconstruction}{cutting} a Morse map
along an arc on a page);
Murasugi sums provide an alternative approach to Goda's results,
Cor.~\ref{subadditivity over Murasugi sum--cor}.

In Sections~\ref{freeness} and \ref{annuli}, 
Murasugi sums and our other constructions 
are used to relate the Morse--Novikov number 
to other knot invariants, as follows.
Let $K$ be a knot, $\b\in B_n$ an $n$-string braid.
The free genus $\freegenus{K}$ of $K$ 
is the least genus of a Seifert surface $S$ 
for $K$ for which $\pi_1(\diff{S^3}{S})$ is a free group.  
The braid index $\bi(K)$ is
the least $m$ such that $K$ can be represented
as a closed $m$-string braid.
The $k$-twisted, $\pm$-clasped Whitehead double of $K$
is the knot $\double{K}{k}{\pm}$ bounding the Seifert surface
$\AKn{K}{k}\plumb{} \AKn{O}{\mp 1}$ plumbed along transverse arcs
of an annulus $\AKn{K}{k}$ having Seifert matrix $\Seif{k}$
with $K\sub\Bd\AKn{K}{k}$ and a Hopf annulus $\AKn{O}{\mp 1}$. %
\label{AKn definition} The wrapping genus $\wrapgenus{K}$ 
(resp., layered wrapping genus $\wlapgenus{K}$) of $K$ is 
the least $n$ such that $K$ lies on a Heegaard surface of 
genus $n$ (resp., $K$ is isotopic to a closed $1$-string 
$o_n$-braid, where $o_n$ is the connected sum of $n$ copies 
of $o_1$).  The crossing number $c(K)$ of $K$ is the least 
number of crossings in a knot diagram for $K$.  

\begin{free genus estimate}\label{free genus estimate statement}
\hyperlink{free genus estimate theorem}%
{$\MN(K)\le 4\freegenus{K}$.}
\end{free genus estimate}
\begin{braid index estimate}
\hyperlink{braid index estimate theorem}%
{$\MN(\double{K}{m}{\pm})\leq 4\bi(K)-2$.}
\end{braid index estimate}
\begin{wrapping genera estimate}
\hyperlink{wrapping genera estimate theorem}%
{$\MN(\double{K}{m}{\pm}\leq 2(\wlapgenus{K}+1)\leq 2(\wrapgenus{K}+1)$.}
\end{wrapping genera estimate}
\begin{crossing number estimate}\label{crossing number estimate statement}
\hyperlink{crossing number estimate theorem}%
{$\MN(\double{K}{m}{\pm}\leq 2(c(K)+2)$.}
\end{crossing number estimate}

In some cases, an upper bound deduced from one of 
these estimates coincides with the lower bound from \cite{P-R-W}, 
so $\MN(L)$ is known precisely.  More often,
unfortunately, a large gap remains: the strongest 
inequality accessible for any knot $K$ using the results 
of \cite{P-R-W} is $\MN(K)\ge 2g(K)$, and we are not aware 
of any technique which could be used to show that 
$\MN(K)>2g(K)$ for some $K$.

\begin{question} Does there exist a knot $K$ with
$\MN(K)>2g(K)$? 
\end{question}

This paper supercedes, and considerably extends, the
second author's preprint \cite{Rudolph:Msums} (in particular,
the proof of the {Free Genus Estimate} in \cite{Rudolph:Msums} 
was inadequate, and \cite{Goda:suture-sum,Goda:handle-number}
had been overlooked).  
Both authors thank Hiroshi Goda, 
Andrei Pajitnov, 
and Claude Weber
for helpful conversations and communications, 
and Walter Neumann for comments 
on a draft of this paper.
The Section de Math\'e\-ma\-tiques of the University of Geneva 
provided extensive hospitality during much of this research.

\section{Preliminaries and generalities}\label{preliminaries}

The symbol $\qedsymbol$ signals either 
the end or the omission 
of a proof. 
The notations $A\isdefinedas B$ and $B\whichdefines A$ 
both define $A$ to mean $B$.  
Terms being defined are set in \bydeffontname~type;
definitions labelled as such are either less standard 
or of greater (local) significance than definitions made 
in passing.

Spaces, maps, etc., are smooth ($\mathscr C^\infty$) 
unless otherwise stated.
The set of critical points of a map $f\from M\to N$
is denoted $\crit f$; for $x\in\crit f$, let $\ind{f}{x}$ 
denote the index of $f$ at $x$.
Manifolds may have boundary, but corners only if so noted.
A \bydef{closed manifold} is one which is compact 
and has empty boundary.
Manifolds are (not only orientable, but) oriented
unless otherwise noted; in particular, $\R$, $\C^n$, 
$ D^{2n}\isdefinedas \{(z_1,\dots,z_n)\in\C^n
\Suchthat
|z_1|^2+\dots+|z_n|^2\le 1\}$,
and $S^{2n-1}\isdefinedas \Bd D^{2n}$ are equipped with standard 
orientations, as is $S^2$ when it is identified with the Riemann sphere 
$\Cext\isdefinedas \C\cup\{\infty\}$.  The manifold $M$ with its
orientation reversed is denoted $-M$.  The interior (resp., boundary)
of $M$ is denoted by $\Int M$ (resp., $\Bd M$).

For suitable $Q\sub M$, let $\Nb{M}{Q}$ denote a closed 
regular neighborhood of $Q$ in $(M,\Bd M)$,
and let $\Ext{M}{Q}$ denote $\diff{M}{\Int\Nb{M}{Q}}$,
the \bydef{exterior} of $Q$.
A submanifold $Q\sub M$ is \bydef{proper} if 
$\Bd Q=Q\cap\Bd M$.
If $Q$ is a codimension--$2$ submanifold of $M$ 
with trivial normal bundle, then a trivialization 
$\t\from Q\times D^2 \to \Nb{M}Q$ 
is \bydef{adapted to} a map $f\from \diff{M}{Q}\to S^1$ if
$\t(\x,0)=\x$ and $f(\t(\x,z)) = z/|z|$ 
for $z\ne 0$.

An \bydef{arc} is a manifold diffeomorphic to $\clint{0}{1}$.
A \bydef{surface} is a compact $2$--manifold.
A \bydef{link} $L$ is a non-empty closed 
$1$--submanifold of $S^3$; a \bydef{knot} is a connected link.
A \bydef{spanning surface} for a link $L$ is a surface $S\sub S^3$ 
with $\Bd S=L$; a \bydef{Seifert surface} for $L$ is a spanning
surface for $L$ without closed components (every $L$ has a Seifert 
surface).  

If $K$ is a knot and $k\in\N$, then a \bydef{$k$-twisted annulus
of type $K$} is any annulus $\AKn{K}{k}\sub S^3$ such that 
$K\sub\Bd\AKn{K}{k}$ and the linking number in $S^3$ of the 
$1$-cycles $K$ and $\diff{\Bd\AKn{K}{k}}{K}$ is $-k$.  
If $S\sub S^3$ is a surface and $K\sub S$ is a knot, then 
the \bydef{$S$-framing of $K$} is the integer $k$ such that
$\Nb{S}{K}=\AKn{K}{k}$.

A \bydef{handlebody of genus $g$} is a boundary-connected sum 
$(S^1\times D^2)_1\bdconnsum\dotsm\bdconnsum(S^1\times D^2)_g%
\linebreak[1]\whichdefines\hb{g}$ of $g\ge 0$ solid tori. 
A handlebody $\hb{g}\sub S^3$ is \bydef{Heegaard} if 
$\diff{S^3}{\Int{\hb{g}}}$ is a handlebody; 
a \bydef{Heegaard surface} is the boundary
of a Heegaard handlebody.  
A \bydef{genus-$g$ Heegaard splitting of $S^3$}
is a pair $(\hb{g}_1,\hb{g}_2)$ where $\hb{g}_1$ is Heegaard
(so $\hb{g}_2$ and $\Bd\hb{g}_1=\Bd\hb{g}_2$ are Heegaard as well).  
According to Waldhausen \cite{Waldhausen1968}, up to isotopy 
there is only one genus-$g$ Heegaard handlebody or surface in,
or splitting of, $S^3$.

Let $L\sub S^3$ be a link.  The image 
of the fundamental class $[S] \in H_2(S,L;\Z)$ 
in $H_2(S^3,L;\Z) \iso H^1(\diff{S^3}{L};\Z)$ 
is independent of the choice of spanning surface 
$S$ for $L$; let $\xi_L\in H^1(\diff{S^3}{L};\Z)\iso %
\pi_0(\operatorname{Map}(\diff{S^3}{L},S^1))$ 
correspond to $[S]$.
Call $f\from \diff{S^3}{L} \to S^1$ \bydef{simple}
if $\xi_L\iso [f]\in \pi_0(\operatorname{Map}(\diff{S^3}{L},S^1))$, 
and \bydef{Morse} if it is smooth and has no degenerate
critical points.
The \bydef{Morse--Novikov number} of $L$, written 
$\MN(L)$, is the least $n$ such that 
some simple Morse map $f\from \diff{S^3}{L} \to S^1$
has $n$ critical points.

\begin{definitions}\label{properties of Morse maps}
Let $f\from \diff{S^3}{L} \to S^1$ be a simple Morse map.
The \bydef{binding} of $f$ is the link $L$.
A \bydef{page} of $f$ is any 
$\page{f}{\theta}\isdefinedas L\cup f^{-1}(\exp(\imunit\theta))$ 
for $\exp(\imunit\theta)\in S^1$; 
the page $\page{f}{\theta}$ is \bydef{smooth} if
$\exp(\imunit\theta)\in\diff{S^1}{\crit{f}}$,
\bydef{singular} if $\exp(\imunit\theta)\in\crit{f}$.
Say that $f$ is:
\begin{inparaenum}[(a)]
\item 
\bydef{boundary-regular} if $f$ has an adapted
trivialization $\t\from L\times D^2 \to \Nb{S^3}L$; 
\item \bydef{moderate} if $\ind{f}{x}\in\{1,2\}$
for all $x\in\crit f$;
\item \bydef{self-indexed} if $f\restr\crit f$ factors as 
$x\mapsto\ind{f}{x}$ followed by an injection;
\item 
\bydef{minimal} if $\card{\crit f}\le\card{\crit g}$ for all 
simple Morse $g\from \diff{S^3}{L} \linebreak[2]\to S^1$;
\item
\bydef{boundary-connected} if every page of $f$ has trivial
second homology (equivalently, if no page of $f$ contains
a non-empty closed surface);
\item
\bydef{connected} if every page of $f$ is connected.
\end{inparaenum}
\end{definitions}

\begin{prop}\label{implications among properties of Morse maps}
Let $f\from \diff{S^3}{L} \to S^1$ be a simple 
Morse map.  
\thmlistfix
\begin{enumerate}
\item\label{finitely many cps is boundary-regular} 
If $\card{\crit f}<\infty$, then up to proper isotopy
$f$ is boundary-regular. 
\item\label{boundary-regular has finitely many critical points}  
If $f$ is boundary-regular, then: 
	\begin{inparaenum}
	\item\label{even number of critical points}
	$\card{\crit f}$ is finite and even;
	\item
    every smooth page of $f$ is a spanning surface for $L$; and
	\item
	$f$ is boundary-connected if and only if every smooth page of $f$
    is a Seifert surface for $L$.
    \end{inparaenum}
\item
If either 
    \begin{inparaenum}
    \item\label{minimal is moderate} 
    $f$ is minimal or 
    \item\label{boundary-connected is moderate}
    $f$ is boundary-connected,
    \end{inparaenum}
then $f$ is moderate.
\item\label{moderate is self-indexed} 
If $f$ is moderate, then up to isotopy $f$ is moderate 
and self-indexed; if also $\card{\crit f}<\infty$, then
the isotopy may be taken to be proper.
\item\label{connected implies boundary-connected}
If $f$ is connected, then $f$ is boundary-connected.  
\item\label{boundary-connected implies connected, sometimes}
If $f$ is boundary-connected and every Seifert surface for $L$ 
is connected \textup(e.g., if $L$ is a knot\textup), then
$f$ is connected.
\item\label{self-indexed} Let $f$ be 
boundary-regular,
moderate, 
and self-indexed.
	\begin{inparaenum}
	\item\label{fibration}
	If $\crit f=\emptyset$, then:
		\begin{inparaenum}
		\item
		$f$ is connected;
		\item
		$f$ is a fibration over $S^1$; and 
		\item\label{fibration:isotopic}
		any two of the pages of $f$ are isotopic 
        \textup{(}rel.\ $L$\textup{)}.
		\end{inparaenum}
	\item\label{non-fibration} 
	If $\crit f\ne\emptyset$, then:
		\begin{inparaenum}
		\item
		half the critical points of $f$ are of index $1$
		and half of index $2$;
		\item\label{two fibrations}
		$f\restr f^{-1}(\diff{S^1}{f(\crit f)})$
		is a trivial fibration over each of the two components of
		$\diff{S^1}{f(\crit f)}$; 
        \item\label{two isotopy classes}
        the smooth pages of $f$ fall into two isotopy classes 
        \textup{(}rel. $L$\textup{)}; and
		\item\label{euler characteristics}
        if $\page{f}{\theta_1}$ and $\page{f}{\theta_2}$ 
        belong to these two isotopy classes, then 
        $|\chi(\page{f}{\theta_1})-\chi(\page{f}{\theta_2})|%
        =\card{\crit f}$.
        \end{inparaenum}
    \end{inparaenum}
\end{enumerate}
\end{prop}
\begin{proof} 
Straightforward.  (For \eqref{minimal is moderate}
and \eqref{non-fibration}, see \cite{P-R-W}.  A slightly
more precise statement of \eqref{euler characteristics} appears 
in Cor.~\ref{from small to large--general}.)
\end{proof}

In case \eqref{fibration}, $L$ is (as usual) called 
a \bydef{fibered link} and $\page{f}{\theta}$ is called 
a \bydef{fiber surface of $f$}.  (Also as usual, 
any Seifert surface for $L$ isotopic to a page of $f$ 
is called a \bydef{fiber surface for $L$}.)
In case \eqref{non-fibration}, 
any smooth page of smaller (resp., larger) Euler 
characteristic will be called a \label{large--def}\bydef{large} 
(resp., \bydef{small}) page, spanning surface, or Seifert surface, as 
the case may be, of $f$ (note: not ``of $L$'').  Any fiber surface 
of a fibration $f$ may be called either large or small, 
as suits convenience.

\begin{convention}\label{Morse map convention}
Henceforth, all Morse maps are boundary-regular and simple.
\end{convention}

\begin{prop}\label{sufficient for bd-connectedness}
If every Seifert surface for $L$ is connected, then 
every minimal Morse map $f\from \diff{S^3}{L} \to S^1$ is 
boundary-connected \textup(and so, by 
Prop.~\textup{\ref{implications among properties of Morse maps}%
\eqref{boundary-connected implies connected, sometimes}}, connected\textup).
\end{prop}
\begin{proof} 
Given any link $L$, and a Morse map 
$f\from \diff{S^3}{L} \to S^1$ which is not 
boundary-connected, there is a spanning surface 
$\page{f}{\theta_0}$ such that the union $S'$ of all 
closed components of $\page{f}{\theta_0}$ is non-empty.  
If the Seifert surface $S''\isdefinedas \diff{\page{f}{\theta_0}}{S'}$ 
is connected, then by Alexander duality 
$H_2(\diff{S^3}{S''};\Z)\iso \reduced{H}_0(S'';\Z)=\{0\}$,
so by a standard argument there is a 
compact $3$--submanifold $M\sub \diff{S^3}{S''}$
such that $\Bd M\ne\emptyset$ is a component of $S'$.
Every $1$--cycle in $M$, being disjoint from $S''$, 
has linking number $0$ with $L$.  It follows
that the restriction $f\restr M\from M\to S^1$
has a continuous lift 
through $\R\to S^1\suchthat\theta\mapsto\exp(\imunit\theta)$
to $\liftof{f}\from M\to \R$;
$\liftof{f}$ is a Morse function rel.\ $S'$, and has the
same critical points, with the same indices, as $f\restr M$.  
Since $M$ is compact, $\liftof{f}$ has (global) extrema in $\Int M$, 
which are local extrema of $f\restr\Int M$ and thus of $f$, 
so $f$ is not moderate.  By 
Prop.~\ref{implications among properties of Morse maps}%
\eqref{minimal is moderate}, $f$ is not minimal.
\end{proof}

For many links $L$ (e.g., knots, fibered links)
the hypothesis of Prop.~\ref{sufficient for bd-connectedness}
is satisfied.  However, for many other links (e.g., split links) 
it fails; for at least some such links, the conclusion
of Prop.~\ref{sufficient for bd-connectedness} also fails
(see Example~\ref{non-boundary-connected Morse map for U}).

\begin{questions}
\begin{inparaenum}
\item
Does there exist a link $L$ for which 
the conclusion of Prop.~\ref{sufficient for bd-connectedness} 
holds although the hypothesis fails?
\item
Does there exist a link $L$ for which no minimal Morse 
function is boundary-connected?
\end{inparaenum}
\end{questions}

\section{Morse maps from Milnor maps}\label{Morse maps from Milnor maps}
The first explicit Morse maps (in fact, fibrations)
for an infinite class of links were given by 
Milnor's celebrated Fibration Theorem \cite{Milnor:singular-points}, 
where they appear as (the instances for $n=2$ of) what 
are now called the ``Milnor maps'' associated
to singular points of complex analytic functions $\C^n\to\C$.  
For present and future purposes, it is useful to extend 
somewhat the framework in which Milnor studied these maps.
Given a non-constant meromorphic function 
$F\from M\meromapto \Cext$
on a complex manifold $M$, let $\divisor{F}$ be the (possibly singular)
complex hypersurface which is the closure in $M$ of 
$F^{-1}(0)\cup F^{-1}(\infty)$. 

\begin{definitions}
The \bydef{argument} of $F$ is 
$\arg(F)\isdefinedas F/|F|\from \diff{M}{\divisor{F}} \to S^1$.
For $M=\C^n$, the \bydef{Milnor map of $F$} is
$\Mmap{F}\isdefinedas\arg(F)\restr(\diff{S^{2n-1}}{\divisor{F}})$; 
for $r>0$, the \bydef{Milnor map of $F$ at radius $r$}, 
denoted by $\Mmapr{F}{r}$, is the Milnor map of 
$(z_1,\dots,z_n)\mapsto F(z_1/r,\dots,z_n/r)$.
\end{definitions}

\begin{lem}\label{Milnor's lemma}
Let $F\from \C^n\meromapto\Cext$ be meromorphic
and not constant.
A necessary and sufficient condition for 
$(z_1,\dots,z_n)\in \diff{S^{2n-1}}{\divisor{F}}$
to be a critical point of $\Mmap{F}$
is that the complex vectors 
\[
\smash[t]{%
(\conjug{z}_1,\dots,\conjug{z}_n), \medspace
\frac{1}{\imunit\thinspace F(z_1,\dots,z_n)}%
\bigl(\frac{\partial F}{\partial z_1}(z_1,\dots,z_n),\dots,
\frac{\partial F}{\partial z_1}(z_1,\dots,z_n)\bigr)
}
\]
be linearly dependent over $\R$.
\end{lem}
\begin{proof}
For holomorphic $F$, this is \cite[Lemma~4.1]{Milnor:singular-points},
and Milnor's proof there applies equally well to meromorphic $F$. 
\end{proof}

\begin{construction}[Milnor maps]
\label{Milnor maps construction}
\hypertarget{Milnor maps hyperconstruction}{}
Let $F\from \C^2\meromapto\Cext$ be meromorphic, not constant, 
and suppose each irreducible analytic component of $\divisor{F}$ 
has multiplicity $1$ (i.e., $F$ has no repeated factors in the 
algebra of meromorphic functions $\C^2\meromapto\Cext$).  
If we let 
\[
m(F)\isdefinedas\inf\{|z|^2+|w|^2\Suchthat (z,w)\in\divisor{F}\}=
\sup\{r\Suchthat rS^3\cap\divisor{F}=\emptyset\},
\]
then the reasoning in \cite{Milnor:singular-points} shows 
that there is a set $X(F)\sub \opint{m(F)}{\infty}$ of radii $r$,
finite in case $F$ is rational, 
and discrete in $\clopint{m(F)}{\infty}$
in any case, such that:
\begin{inparaenum}[(a)]
\item\label{almost all L(F,r) are links}
if $r\in \diff{\opint{m(F)}{\infty}}{X(F)}$, 
then $\divisor{F}$ intersects $rS^3$ transversally, 
so that $L(F,r)\isdefinedas (1/r)(\divisor{F}\cap rS^3)$ 
is a link in $S^3$;
\item\label{isotopic Milnor links}
if $r$ and $r'$ are in the same component 
of $\diff{\opint{m(F)}{\infty}}{X(F)}$ then 
$L(F,r)$ and $L(F,r')$ are isotopic; and 
\item\label{Milnor maps are simple and adaptable} 
if $m(F)<r\not\in X(F)$, then $\Mmapr{F}{r}$ is simple, 
and there is a trivialization 
of $\Nb{S^3}{L(F,r)}$ which is adapted to $\Mmapr{F}{r}$---so 
that, if also $\Mmapr{F}{r}$ has no degenerate critical points, 
then $\Mmapr{F}{r}$ is a Morse map (in the sense of the 
convention on p.~\pageref{Morse map convention}).
\end{inparaenum}
In practice, Lemma~\ref{Milnor's lemma}
makes it easy to locate the critical points of a Milnor map
and check them for non-degeneracy.

Several special cases of this construction have special names.

If $F$ is holomorphic and $(0,0)\in\divisor{F}=F^{-1}(0)$ 
(so $m(F)=0$), then for every $r\in\opint{0}{\inf X(F)}$ the
link $L(F,r)$ (known as the \bydef{link of the singularity} 
of $F$ at $(0,0)$) is fibered, and $\Mmapr{F}{r}$ is a 
fibration (known as the \bydef{Milnor fibration} of $F$).  
Up to isotopy, the link of the singularity and the Milnor fibration 
are independent of $r<\inf X(F)$.

If $F$ is meromorphic, $(0,0)\in\divisor{F}$ (so $m(F)=0$), 
and neither $F$ nor $1/F$ is holomorphic at $(0,0)$ 
(that is, $(0,0)$ is a \bydef{point of indeterminacy} of 
$F$: it is in both the closure of $F^{-1}(0)$ and the 
closure of $F^{-1}(\infty)$), then we will call $L(F,r)$ 
the \bydef{link of indeterminacy} of $F$ at $(0,0)$.  
Up to isotopy, the link of indeterminacy is independent 
of $r<\inf X(F)$.

If $F$ is a polynomial (i.e., both holomorphic and rational),
then up to isotopy the link $L(F,r)$ is independent of 
$r>\sup X(F)$; it is called the \bydef{link at infinity} 
of $F$ (cf.\ \cite{embeddings}) and denoted $L(F,\infty)$.
(Warning: in \cite{Neumann-Rudolph}, \cite{Neumann}, 
and elsewhere, the phrase ``link at infinity of $F$''
denotes what we will call the \bydef{generic link at 
infinity} of $F$, that is, the intersection of a generic 
fiber of $F$ with any sufficiently large $3$--sphere;
that link can differ substantially from the link at infinity
of $F$ defined here.)

If $F$ is rational and neither a polynomial nor the 
reciprocal of a polynomial, then up to isotopy $L(F,r)$ 
is independent of $r>\sup X(F)$; we will call it the 
\bydef{link of indeterminacy at infinity} of $F$ and 
denote it $L(F,\infty)$.
\end{construction}
Several of the following examples of Milnor maps
will be used in later sections.

\begin{example}\label{o{p,q} as a Milnor fibration}
Given $p,q\ge 0$ with $p+q=1$ in case $pq=0$,
let $F_{p,q}\from \C^2\to\C\suchthat (z,w)\mapsto pz^p+qw^q$; 
$X(F_{p,q})=\emptyset$ because $F_{p,q}$ is weighted-homogeneous.
The link of the singularity of $F_{0,1}^{-1}(0)$ at $(0,0)$
(of course $F_{0,1}^{-1}(0)$ is not in fact singular at $(0,0)$),
namely, $\{(z,w)\in S^3\Suchthat w=0\}$, is denoted $O$, 
and the Milnor fibration of $F_{0,1}=\pr_2$ will 
be denoted $o$.  More generally, the link of the 
singularity of $F_{p,q}^{-1}(0)$ at $(0,0)$ is denoted 
$O\{p,q\}$, and the Milnor fibration of $F_{p,q}$ 
will be denoted $o\{p,q\}$.  
\end{example}

\begin{definitions}\label{unknots definition}
A knot isotopic to $O$ is called an \bydef{unknot};
in particular $O$ itself is called the \bydef{horizontal}
unknot, and $O'\isdefinedas \{(z,w)\in S^3\Suchthat z=0\}$
is called the \bydef{vertical} unknot.
A link isotopic to $O\{p,q\}$ (resp., to the mirror 
image of $O\{p,q\}$) is called a \bydef{torus link} 
of type $(p,q)$ (resp., type $(p,-q)$, or equivalently
$(-p,q)$).  A torus link of type $(2,2)$ (resp., $(2,-2)$)
is a \bydef{positive} (resp., \bydef{negative}) \bydef{Hopf link},
and its fiber surface $\AKn{O}{-1}$ is a \bydef{positive}
(resp., \bydef{negative}) \bydef{Hopf annulus}.
\end{definitions}

A link of indeterminacy is never a knot, and is always
obtained from a (multicomponent) link of a singularity 
by reversing the orientation of some, but not all, 
components; a link of indeterminacy is never isotopic 
to the link of a singularity.  Recent work of Pichon 
\cite{Pichon} confirms an empirical observation (about 
certain real polynomial maps $\R^4\to\R^2$) recorded offhandedly
at the end of \cite[Example 4.7]{icp1}, and immediately implies 
that the link of indeterminacy of $F$ is a fibered link if and only if 
$\Mmapr{F}{r}$ is a fibration for all sufficiently small $r>0$.

\begin{example}\label{fibered and unfiberable links of indeterminacy}
Calculations using Lemma~\ref{Milnor's lemma} show that the 
link of indeterminacy of $(z^2+w^3)/(z^3+w^2)$ at $(0,0)$
is fibered by the Milnor map at small radius.  
The link of indeterminacy of $z/w$ at $(0,0)$ (and at infinity)
is a negative Hopf link and is fibered by its Milnor map. 
The Milnor map of $(z^2+w^2)/(z^2-w^2)$ at any radius $r>0$ 
is Morse but not a fibration; the link of indeterminacy 
at $(0,0)$ (and at infinity) is not fiberable, since it has a 
disconnected Seifert surface (two disjoint annuli), whereas for 
homological reasons every Seifert surface of a fibered link is connected.
\end{example}

Like links of indeterminacy, many links at infinity are 
fibered (for instance, if $L(F,\infty)$ is a knot, then
it is fibered \cite{Neumann-Rudolph}) but some are 
non-fiberable; in the fibered case, the Milnor map 
at sufficiently large radius is a fibration, and in 
every case the Milnor map is Morse provided that its 
critical points are non-degenerate.
There is a very limited overlap between the classes
of links at infinity and links of singularities (for
instance, since $X(F_{p,q})=\emptyset$, the link at
infinity of $F_{p,q}$ is $O\{p,q\}$, and---as it 
happens---any knot which is both a link at infinity and
the link of a singularity is a torus knot of type $(p,q)$ 
for positive, relatively prime $p$ and $q$), \cite{embeddings,Neumann}.
The polynomials with fibered links at infinity have been
variously characterized by several authors (see Bodin \cite{Bodin}
and references therein).

\begin{example}\label{fibered and non-fiberable links at infinity}
Both $F\Suchthat(z,w)\mapsto (w^2-z)^2-z^5-4z^4 w$ and
$G\Suchthat(z,w)\mapsto z(zw+1)$ 
have links at infinity which are fibered by 
their Milnor maps at sufficiently large radius.  
The knot $L(F,\infty)$ is not the link of a singularity:
$L(F,\infty)$ is a non-trivial cable on a non-trivial torus 
knot (as is readily seen from the parametrization 
$\zeta\mapsto(\zeta^4,\zeta^5+\zeta^6)$ of $F^{-1}(0)$),
and therefore not itself a torus knot.  
The generic link at infinity of $G$ is the link at infinity of 
$(z,w)\mapsto(z+t)((z+t)w+1)-t = z^2 w + 2tzw + z + t^2w$ for $t\ne 0$,
and is not fiberable (see \cite[corrigendum]{Neumann-Rudolph}).
\end{example}

Like a link of indeterminacy, a link of indeterminacy at infinity
is never a knot, and is always obtained from a (multicomponent) 
link at infinity by reversing the orientation of some, but not all, 
components.  A link of indeterminacy at infinity can be isotopic 
to a link at infinity.  Some links of indeterminacy at infinity 
are fibered, others are non-fiberable.

\begin{question}
If the link of indeterminacy at infinity of $F$ is fibered,
is $\Mmapr{F}{r}$ is a fibration for all sufficiently large $r$?
(Can the methods of \cite{Pichon} be adapted to this context?)
\end{question}

\begin{example}\label{fibered loi@infinity}
Let $G_0\from\C^2\meromapto\Cext\suchthat (z,w)\mapsto zw/(4z-1)$.
Easily, $X(G_0)=\{1/4\}$.  Calculation with Lemma~\ref{Milnor's lemma}
shows that, for every $r>0$, $\Mmapr{G_0}{r}$ has no critical points.  
Thus the link of indeterminacy at infinity $L(G_0,\infty)=L(G_0,1)$ 
is fibered by $\Mmap{G_0}$.  
(Alternatively, though less explicitly, it is easy to recognize 
that $L(G_0,\infty)$ is isotopic to the connected sum of a
positive and a negative Hopf link; and the connected sum of fibered 
links is well known to be fibered---a proof from the point of view of 
Morse maps is given in \cite{P-R-W}.)  More generally, if $k\in\Z$ 
and $G_k\from\C^2\meromapto\Cext\suchthat (z,w)\mapsto zw/(4z-w^k)$,
then $L(G_k,\infty)=L(G_k,1)$ is fibered by $\Mmap{G_k}$.  
(Notice that $L(G_1,\infty)$ is isotopic to the link at infinity 
$L(G,\infty)$ in Example~\ref{fibered and non-fiberable links 
at infinity}.)
\end{example}

\begin{example}\label{non-solvable} 
\begin{figure}
\centering
\includegraphics[width=0.25\textwidth]{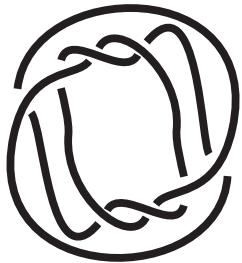}
\caption{\label{non-solvable:figure} Can this fibered knot 
be fibered by a Milnor map?}
\end{figure}
It can be shown (for instance, by using the techniques 
of \cite{algfuns}) that, for all sufficiently small $\e>0$,
if $F(z,w)=1-z^2+3z^6+(\e w)^3-3\e w$, then $L(F,1)$ is 
the $8$-crossing knot pictured in Fig.~\ref{non-solvable:figure}.  
This knot is not the link of a singularity nor a link at infinity.
Although it is fibered (being a closed positive braid, 
\cite{Stallings}), calculations using Lemma~\ref{Milnor's lemma} 
show that $\Mmap{F}$ has at least $2$ critical points.  We do 
not know if there exists $G$ such that some link $L(G,r)$ is 
isotopic to $L(F,1)$ and $\Mmapr{G}{r}$ is a fibration.
\end{example}

\begin{example}\label{loi@infinity which is l@infinity}
If $F\from\C^2\to\C\suchthat (z,w)\mapsto z(2z-1)$ and
$G\from\C^2\meromapto\Cext\suchthat (z,w)\mapsto z^{-1}(2z-1)=2-z^{-1}$,
then $L(G,\infty)$ is obtained from $L(F,\infty)$ by reversing
the orientation of one component; since both links are evidently
split links of two unknotted components, they are isotopic.  
For $r\ge 1/4$, and in particular for $r>\sup X(F)=1/2$, 
$\Mmapr{F}{r}$ has a $1$--sphere of degenerate critical points,
so is not Morse.  By contrast, $\Mmapr{G}{r}$ is Morse for 
$r\ge 1/2=\sup X(G)$, with one critical point of index $1$ 
and one of index $2$.
Let $U\isdefinedas L(G,1)=L(G,\infty)$,
$u\isdefinedas\Mmap{G}$.
(The situation is pictured in Fig.~\ref{uround}.)  
\begin{figure}
\centering
\includegraphics[width=0.5\textwidth]{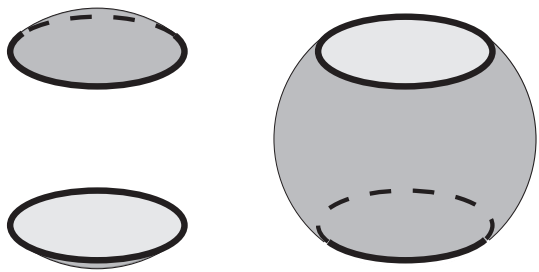}
\caption{\label{uround} A small and a large 
Seifert surface of $u$.}
\end{figure}
\end{example}

\begin{example}\label{o1 as a Milnor map}
A Milnor map which is a Morse map need not be minimal.
Let $F(z,w)=4w+3(w^2+z^2)$.  
Calculations using Lemma~\ref{Milnor's lemma} 
show that $X(F)=\{4/3\}$; for $0<r<4/3$, $L(F,r)$ is
isotopic to the link of the singularity of $F$ at $(0,0)$,
an unknot, while for $4/3<r<\infty$, $L(F,r)$ is isotopic
to $O\{2,2\}$, the link at infinity of $F$.  
For $0<r<2/3$ or $4/3<r$, $\Mmapr{F}{r}$ is a fibration;
for $r=2/3$, $\crit{\Mmapr{F}{r}}=%
\{(z,w)\Suchthat \Re(z)=0,\Re(w)=-1/3, |z|^2+|w|^2=4/9\}$
is a circle of degenerate critical points; 
for $2/3<r<4/3$, $\Mmapr{F}{r}$ is Morse, 
and $\crit{\Mmapr{F}{r}}=\{(0,w)\Suchthat |w|=r, |w+1|=1/3\}$
consists of two points, one of index $1$ and one of index $2$.
In particular, $\Mmap{F}$ is a Morse map for 
the unknot $L(F,1)$.  Let $o_1\from\diff{S^3}{O}\to S^1$ 
be the Morse map for $O$ onto which $\Mmap{F}$ is 
carried by an isotopy carrying $L(F,1)$ onto $O$.
\begin{figure}
\centering
\includegraphics[width=0.5\textwidth]{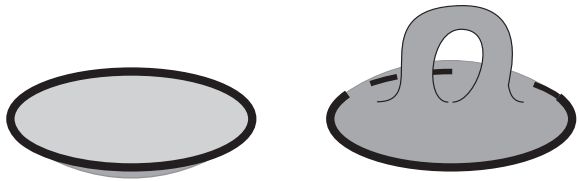}
\caption{\label{o1} A small and a large Seifert surface of $o_1$.}
\end{figure}
A small and a large Seifert surface of $o_1$ are pictured in Fig.~\ref{o1}.
\end{example}
The non-fibrations $u$ and $o_1$, trivial though they be, 
are ingredients of fundamental importance throughout the
following sections.

\begin{prop}\label{MN(U)=2}
$\MN(U)=2$ and $u$ is minimal.
\end{prop}

\begin{proof}
Since $U$ has a disconnected Seifert surface, it is not fibered.
\end{proof}

\begin{example}\label{non-boundary-connected Morse map for U}
If $\MN(L)=0$, then any two minimal Morse maps for $L$ are
isotopic: this simply rephrases the well-known fact that, 
up to isotopy, a fibered link has a unique fibration.  
If $\MN(L)>0$, then there may be more than one isotopy class
of minimal Morse map.  We illustrate this by constructing 
a minimal Morse map for $U$ which is not isotopic to $u$.

First, construct two Morse maps---non-moderate, 
and barely non-minimal---for $O$, as follows.
Let $B^3\sub\diff{S^3}{O}$ be a ball with 
$1\in o(\Int B^3)$ and $-1\ne o(B^3)$, and such that  
(a continuous branch of) $-\imunit\log(o|\Bd B^3)$ is a Morse
function on $\Bd M$ with exactly two critical points.  
By a standard argument (see \cite{Milnor:Morse}), 
for ${\ext}\in\{{\min},{\max}\}$ there is a smooth homotopy
(supported off a neighborhood of $\Bd B^3$) 
from $-\imunit\log(o|B^3)$ to a function 
$g_{\ext}\from B^3\to\opint{-\pi/2}{\pi/2}$ such that:
\begin{inparaenum}[(a)]
\item
$g_{\ext}$ has exactly two critical points in $\Int B^3$, 
both non-degenerate, having indices $0$ and $1$ in the 
case of $g_{\min}$ and indices $2$ and $3$ in the case 
of $g_{\max}$;
\item
$0$ is a regular value of $g_{\ext}$;
\item
the critical values of $g_{\ext}$ are of opposite sign.
\end{inparaenum}
Then there is a smooth homotopy 
(supported in $\Int B^3$) from
$o$ to a non-moderate Morse map 
$o_{\ext}\from\diff{S^3}{O}\to S^1$ 
with $o_{\ext}\restr B^3=\exp(\imunit g_{\ext})$
and $o_{\ext}(B^3)=o(B^3)$.  
Clearly $\page{o_{\ext}}{-1}$ is a $2$--disk 
and $\page{o_{\ext}}{1}$ is the disjoint 
union of a $2$--disk and a $2$--sphere $S^2_{\ext}$.
Let $B^3_{\ext}\sub S^3$ be the $3$--ball with 
$\Bd B^3_{\ext}=S^2_{\ext}$ and 
$\page{o_{\ext}}{-1}\sub\Int B^3_{\ext}$.
The identification space 
$\Sigma\isdefinedas (B^3_{\min}\sqcup B^3_{\max}){/}{\equiv}$,
where $\equiv$ identifies $S^2_{\min}$ to $S^2_{\max}$ by a 
diffeomorphism, is a piecewise-smooth $3$--sphere which
is easily given a smooth structure such that
\begin{multline*}
\upsilon\isdefinedas 
\left(o_{\min}\restr(\diff{B^3_{\min}}{O}) \sqcup %
o_{\max}\restr(\diff{B^3_{\max}}{O})\right) {/}{\equiv} \\
\from \diff{\Sigma}{(\Bd \page{o_{\min}}{-1}\cup %
\Bd \page{o_{\max}}{-1})}\to S^1
\end{multline*}
has exactly two critical points, both nondegenerate,
of indices $1$ and $2$.  (The same construction 
one dimension lower is pictured in Fig.~\ref{nonbconu}.) %
\begin{figure}
\centering
\includegraphics[width=.9\textwidth]{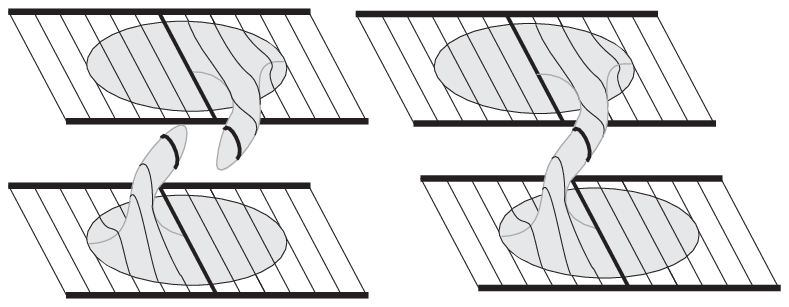}
\caption{\label{nonbconu} 
Top left: a Morse modification,
with a local minimum and a saddlepoint, of \protect$-\imunit$ times
the logarithm of the Milnor fibration of the identity map of
\protect$\diff{\Cext}{\{0,\infty\}}$.  
Bottom left: a similar map
with a saddlepoint and a local maximum.  
Left: \protect$-\imunit$ times
the logarithm of a Morse map (with two saddlepoints) 
\protect$\diff{S^2}{\{x_1,x_2,x_3,x_4\}}
\protect\cong
(\diff{\Cext}{\{0,\infty\}})
\connsum
(\diff{\Cext}{\{0,\infty\}})\to S^1$.
}
\end{figure}%
If $\d\from S^3\to\Sigma$ is a diffeomorphism 
with $\d(U)=\Bd \page{o_{\min}}{-1} \cup \Bd \page{o_{\max}}{-1}$,
then $\omega\isdefinedas \upsilon\after\delta$
is a non-boundary-connected minimal Morse map for $U$;
a large Seifert surface of $\omega$ is the union of two
disjoint $2$--disks, and a small spanning surface is 
the disjoint union of a $2$--sphere and two $2$--disks
which are separated by the $2$--sphere.
\end{example}

Links of singularities, links at infinity, links of indeterminacy,
and links of indeterminacy at infinity are all \bydef{graph links}
in $S^3$ as defined by Eisenbud \& Neumann~\cite{Eisenbud-Neumann}.
Graph links are highly atypical of links $L$ which have Milnor
maps $\diff{S^3}{L}\to S^1$.  A large class of such links 
(namely, precisely those such links whose Milnor maps come from 
holomorphic---not merely meromorphic---maps $F\from\C^2\to\C$)
consists of the \bydef{transverse $\C$--links} in the sense of 
\cite{totaltan}; by Boileau \& Orevkov \cite{Boileau-Orevkov}, 
transverse $\C$-links are exactly the same, up to isotopy, as 
\bydef{quasipositive} links (defined and studied in \cite{constqp1}
and its sequels).  A typical quasipositive link---for example,
the knot in Example~\ref{non-solvable}---is very far from 
being a graph link.  

We conclude this section with a series of questions about links 
which have, up to isotopy, a Milnor map which is a minimal Morse 
map---briefly, links which are \bydef{\Mgood}.  
Prop.~\ref{MN(U)=2} asserts that $U$ is \Mgood; 
Milnor's Fibration Theorem \cite{Milnor:singular-points} 
proves that the link of a singularity is \Mgood; 
Pichon's results \cite{Pichon} show that a fibered link 
of indeterminacy is \Mgood.  Although the knot in 
Example~\ref{non-solvable} has a Milnor map which is Morse,
we do not know if it is \Mgood.

\begin{questions}\label{Milnor map questions}
\begin{inparaenum}
\item\label{does qp mean MM?}
Is every transverse $\C$-link \Mgood?  
\item\label{does strongly qp mean MM?}
If \eqref{does qp mean MM?} cannot be answered in the
affirmative (or until it is), is at least every strongly
quasipositive link \Mgood?  (By definition, $L$ is 
strongly quasipositive if and only if $L$ has a quasipositive 
Seifert surface, as defined in \cite{braided-surfaces}; 
by \cite{constqp3}, $L$ is strongly quasipositive if and only if 
$L$ bounds a subsurface of a fiber surface of a torus link of type
of $\{p,q\}$ for some $p,q \ge 1$.  A strongly quasipositive link 
is quasipositive; many quasipositive links are not strongly quasipositive.)
\item\label{does strongly qp fibered mean MM?}
If~\eqref{does strongly qp mean MM?} cannot be answered in the
affirmative (or until it is), is at least every strongly
quasipositive fibered link \Mgood---that is, does every 
strongly quasipositive fibered link have a fibration which,
up to isotopy, is the Milnor map of a holomorphic function?  
(It follows from \cite{Neumann-Rudolph, icp1,constqp3}, 
and deep results of Giroux \cite{Giroux:fibered,Giroux:ICM}
that a fibered link $L$ is strongly quasipositive 
if and only if $L$ is a stable plumbing of positive Hopf links; 
see p.~\pageref{annular plumbing}.)
\item
Is at least the knot in Fig.~\ref{non-solvable:figure},
which is fibered and strongly quasipositive, \Mgood?  
\end{inparaenum}
\end{questions}

\section{Stallings twists of Morse maps}%
\label{twisting}

\begin{figure}
\centering
\includegraphics[width=0.6\textwidth]{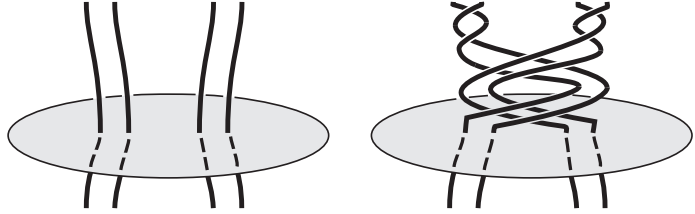}
\caption{\label{Stallings twist:figure} A $2$-disk 
bounded by the unknot $O$, and pieces of links 
$L$ and $\multitwist{L}{O}{1}$ in a neighborhood 
of that $2$-disk.}
\end{figure}

Define $\blowup\from S^3\to S^3$ by 
$\blowup(z,w)=(z\arg(w),w)$ for $w\ne 0$, $\blowup(z,0)=(z,0)$.
Although the bijection $\blowup$ is discontinuous at each point of 
the standard unknot $O$, its restriction
$\blowup\restr(\diff{S^3}{O})\from\diff{S^3}{O}\to\diff{S^3}{O}$
is a diffeomorphism (as is $\blowup\restr O=\id{O}\from O\to O$).
If $X\sub\diff{S^3}{O}$ is a manifold, then the manifold 
$\multitwist{X}{O}{n}\isdefinedas \blowup^n(X)\sub S^3$ 
is diffeomorphic to $X$, but may or may not be isotopic to $X$ 
for $n\ne 0$.  (Fig.~\ref{Stallings twist:figure} illustrates the
relationship among $X$, $O$, and $\multitwist{X}{O}{1}$
when $X=L$ is a link.)
More generally, given any unknot $\g\sub S^3$, and a diffeomorphism
$\d\from S^3\to S^3$ with $\g=\d(O)$, let 
$\blowup_\g\isdefinedas\d^{-1}\after\blowup\after\d$,
and for $X\sub\diff{S^3}{\g}$, let
$\multitwist{X}{\g}{n}\isdefinedas \blowup_\g^n(X)$. 
Up to isotopy, $\multitwist{X}{\g}{n}$ does not depend
on the choice of $\d$.  

Let $M\sub S^3$ and $Q$ be manifolds. 
Say a smooth map $f\from M\to Q$ is
\bydef{flat near $O$} when $O$ has a 
regular neighborhood $\Nb{S^3}{O}$
such that $f\restr(M\cap\Nb{S^3}{O})$ 
factors through $(z,w)\mapsto((1-|w|^2)^{1/2},w)$.  
More generally, given an unknot $\g$ and 
a diffeomorphism $\d\from S^3\to S^3$ with $\g=\d(O)$, 
say that $f$ is \bydef{flat near $\g$} if $f\after\delta^{-1}$
is flat near $O$.  (This does depend on the choice of
$\d$.)  If $f$ is flat near $\g$, then $f$ cannot see 
that $\blowup_\g$ is discontinuous---in fact,
$f\after\blowup_\g\from M\to Q$ is smooth, 
$\crit{f\after\blowup_\g}=\blowup_\g^{-1}(\crit{f})$,
and if $Q=\R$ or $Q=S^1$, then 
$\blowup_\g\restr\crit{f\after\blowup_\g}$ is index-preserving.  

\begin{construction}[Stallings twists] 
\label{Stallings twists construction}
\hypertarget{Stallings twists hyperconstruction}{}
Let $L\sub S^3$ be a link, $f\from\diff{S^3}{L}\to S^1$ 
a Morse map. If 
\begin{inparaenum}[(a)]
\item\label{twist on an unknot}
$\g$ is an unknot,
\item\label{unknot is on one page}
for some $\theta$, $\g\sub\diff{\Int{\page{f}{\theta}}}{\crit{f}}$,
and 
\item\label{unknot has framing 0}
the $\page{f}{\theta}$-framing of $\g$ is $0$,
\end{inparaenum}
then it is easy to check that $f$ is flat near $\g$
for an appropriate choice of $\d$.  Having made 
that choice, call $\multitwist{f}{\g}{n}\isdefinedas %
f\after\d\after\blowup^n\from %
\diff{S^3}{\multitwist{L}{\g}{n}}\to S^1$ 
the \bydef{$n$-fold Stallings twist of $f$ along $\g$}.
\end{construction}

\begin{unproved}{prop}\label{twist proposition}
The $n$-fold Stallings twist of $f$ along $\g$
depends \textup(up to isotopy\textup) 
only on $L$ and $f$ \textup(up to isotopy),
$\g$ \textup(up to isotopy respecting property 
\eqref{unknot is on one page}, allowing $\theta$ to vary\textup), 
and $n\in\Z$, and $\multitwist{f}{\g}{n}$ is a Morse map
for $\multitwist{L}{\g}{n}$ with the same number of 
critical points of each index as the Morse map $f$.
\end{unproved}

\begin{unproved}{cor}\label{twist corollary} 
If $L\sub S^3$ is a link, $f\from\diff{S^3}{L}\to S^1$ 
is a Morse map, and $\g\sub\diff{\Int(\page{f}{\theta})}{\crit{f}}$
is an unknot with $\page{f}{\theta}$-framing $0$,
then, for every $n\in\Z$, 
$\MN(\multitwist{L}{\g}{n})=\card{\crit{f}}\le\MN(L)$.
\end{unproved}

\begin{example}\label{fibered loi@infinity, bis}
In Example~\ref{fibered loi@infinity}, 
the $\page{G_0}{\theta}$-framing of $O\sub L(G_0,1)$ 
is $0$ for all $\theta$.  If $\g\isdefinedas\diff{\AKn{O}{0}}{O}$, 
then $\multitwist{L(G_0,1)}{\g}{n}$ is (isotopic to)
$L(G_n,1)$ and $\multitwist{f}{\g}{n}$ is (isotopic to)
$\Mmap{G_n}$, for all $n\in\Z$.
\end{example}

\begin{prop}
If $L\sub S^3$ is a link, $f\from\diff{S^3}{L}\to S^1$ is a Morse map,
and $\g\sub\diff{\Int(\page{f}{\theta})}{\crit{f}}$, then
for every $n$ there is a Morse map $\multitwist{f}{\g}{n}\to S^1$ 
which has exactly two more critical points than $f$, 
one of index $1$ and one of index $2$.
\end{prop}

\begin{proof}
If necessary, change coordinates on $S^1$ (by rotation, say)
to arrange that $\exp(\imunit\theta)$ is such that
$\page{o_1}{\theta}$ is a large page of the Morse map
$o_1\from\diff{S^3}{O}\to S^1$ described in 
Example~\ref{o1 as a Milnor map}.  Let $f\connsum o_1$ 
be a connected sum of $f$ and $o_1$ (as defined, say, in 
\cite{P-R-W}; or see Construction~\ref{M. sum of Morse maps construction})
along any component of $L$ which belongs 
to the boundary of the component of $\page{f}{\theta}$ that 
contains $\g$; $f\connsum o_1$ is Morse, with exactly two 
more critical points than $f$ (one of index $1$ and
one of index $2$), and the page $\page{f\connsum o_1}{\theta}$ 
may be identified with an appropriate boundary-connected sum 
$\page{f}{\theta}\bdconnsum \page{o_1}{\theta}$.
A glance at Fig.~\ref{o1} shows that, for every $k$,
$\page{o_1}{\theta}$ contains an unknot $\g_k$ with 
$\page{o_1}{\theta}$-framing $k$.
If the $\page{f}{\theta}$-framing of $\g$ is $s$,
the unknot $\g\connsum\g_{-s}\sub \page{f\connsum o_1}{\theta}$ 
has $\page{f\connsum o_1}{\theta}$-framing $0$.
Evidently the pairs $(L,\g)$ and $(L\connsum O,\g\connsum\g_{-k})$ 
are isotopic, so the links $\multitwist{L}{\g}{n}$ and 
$\multitwist{L\connsum O}{\g\connsum\g_{-k}}{n}$ are
isotopic for any $n$.  
\end{proof}

\begin{unproved}{cor}\label{arbitrary twist}
For any link $L\sub S^3$, any unknot $\g$ interior to 
a \textup(smooth\textup) page of a minimal Morse map
$\diff{S^3}{L}\to S^1$, and all $n$, 
$\MN(\multitwist{L}{\g}{n})\le\MN(L)+2$.
\end{unproved}

\begin{example}\label{MN Bd A{O,n} at most 2}
Let $L$ be $O\{2,2\}$, the fibered positive Hopf link, 
with fiber surface $\AKn{O}{-1}$, the positive Hopf annulus.
Let $\g\sub\AKn{O}{-1}$ be a core circle.
Evidently $\multitwist{L}{\g}{n}=\Bd\AKn{O}{n-1}$; 
by Cor.~\ref{arbitrary twist}, $\MN(\Bd\AKn{O}{k})\le 2$ for 
all $k$.
\end{example}

\begin{unproved}{cor}\label{MN Bd A{O,n} is 2, except}
$\MN(\Bd\AKn{O}{k})=2$ for $k\ne\pm 1$, and 
$\MN(\Bd\AKn{O}{\pm1})=0$.
\end{unproved}

\begin{hist} 
Stallings \cite{Stallings} observed that, 
if $L$ is a fibered link and $\g$ is an unknot on a fiber 
surface $F$ of $L$ with $F$-framing $0$, 
then $\multitwist{L}{\g}{n}$ is a fibered link for any $n$;
he called the operation in question simply ``twisting''.
In \cite{Harer}, Harer observed that, if $L$ is a fibered 
link and $\g$ is an unknot on a fiber surface $F$ of $L$ 
with $F$-framing $\pm 2$, then $\multitwist{L}{\g}{\mp1}$ 
is a fibered link; he called this operation and the one 
introduced in \cite{Stallings} both ``Stallings twists''.
It is clear that, like Stallings's original twisting, Harer
twists can be extended from fibered links and fibrations 
to arbitrary links and Morse functions, and that 
Prop.~\ref{twist proposition} and Cor.~\ref{twist corollary} 
extend in obvious ways.
\end{hist}

\section{Monodromies and adiexodons for Morse maps}%
\label{monodromies}

Let $L\sub S^3$ be a link, $f\from \diff{S^3}{L} \to S^1$ 
a moderate, self-indexed Morse map.  Assume, changing 
coordinates on $S^1$ if need be, that $\page{f}{0}$ 
(resp., $\page{f}{\pi}$) is a small (resp., large) 
spanning surface of $f$ (as defined on p.~\pageref{large--def}) 
and that, if $x\in\crit{f}$, then 
$f(x)=\exp((-1)^{1+\ind{f}{x}}\imunit\pi/4)$. 
Let $\Hs_\ell\isdefinedas L \cup f^{-1}(\{z\in S^1\Suchthat \Re(z)\le 0\})$,
$\Hs_s\isdefinedas L \cup f^{-1}(\{z\in S^1\Suchthat \Re(z)\ge 0\})$.

\begin{lem}\label{Hs and Hl are smooth} 
$\Hs_s$ and $\Hs_\ell$ are smooth $3$--manifolds.
\end{lem}

\begin{proof} 
Since $\crit{f}\sub \Int\Hs_s$, the (topological) boundary 
of $\Hs_s$ (and of $\Hs_\ell$) is the union 
$\page{f}{\pi/2}\cup \page{f}{-\pi/2}$ 
of two spanning surfaces of $f$ along their common boundary $L$,
and so is a closed topological $2$--manifold.  Although the
manifold structure of $\Bd \Hs_s$ is not \foreign{a priori} 
smooth along $L$ (where there might be a corner), in fact 
the convention that $f$ is boundary-regular ensures that 
$\Bd\Hs_s$ is a closed smooth surface, so $\Hs_s$ and $\Hs_\ell$
are smooth.  (Note that as geometric $2$--cycles, 
$\Bd\Hs_s=\page{f}{\pi/2}-\page{f}{-\pi/2}$ and 
$\Bd\Hs_\ell=-\Bd\Hs_s=-\page{f}{\pi/2}+\page{f}{-\pi/2}$.)
\end{proof}

Let $\Nb{S^3}L$ be a regular neighborhood with a trivialization
$\t\from L\times D^2 \to \Nb{S^3}L$ adapted to $f$, so that 
$S_\ell\isdefinedas \page{f}{\pi}\cap\Ext{S^3}{L}$ 
(resp., $S_s\isdefinedas \page{f}{0}\cap\Ext{S^3}{L}$) is 
isotopic in $\diff{S^3}{L}$ to $\page{f}{\pi}$ (resp., $\page{f}{0}$).
Let $g\isdefinedas \rank H_1(S_\ell;\Z)$,
$g_s\isdefinedas \rank H_1(S_s;\Z)$, so that, by 
Prop.~\ref{implications among properties of Morse maps}.%
\eqref{even number of critical points},
$(g-g_s)/2\whichdefines \nu\ge 0$ is an integer.

\begin{prop}\label{Hs and Hl as cobordisms}
There exists a $3$--manifold with corners,
equipped with a handle decomposition
\begin{equation*}\label{domain of eta-s}
\smash[t]{
S_s\times\clint{-1}{1})\cup \bigcup_{k=1}^{2\nu}\h1k 
}
\tag{$\S$}
\end{equation*}
for which 
the $1$--handle $\h1k$ is attached to $S_s\times\clint{-1}{1}$ 
along $S_s\times \{-1\}$ for $k=1,\dots,\nu$ 
and along $S_s\times \{1\}$ for $k=\nu+1,\dots,2\nu$,
and diffeomorphisms of $3$--manifolds with corners
\begin{gather*}
\smash[t]{\eta_s\from (S_s\times\clint{-1}{1})\cup \bigcup_{k=1}^{2\nu}\h1k %
          \to \Hs_s\cap\Ext{S^3}{L},}
\\
\eta_\ell\from S_\ell\times\clint{-1}{1}\to \Hs_\ell\cap\Ext{S^3}{L}
\end{gather*}
such that
\begin{inparaenum}
\item\label{f on Hs}
$f\restr \eta_s(S_s\times\clint{-\e}{\e}) = 
\exp(\imunit\pi(2+\pr_2\after\eta_s^{-1})/2)$ 
for all sufficiently small $\e>0$, and
\item\label{f on Hl}
$f\restr \Hs_\ell\cap\Ext{S^3}{L} = 
\exp(\imunit\pi(2+\pr_2\after\eta_\ell^{-1})/2)$.
\end{inparaenum}
\end{prop}

\begin{proof} By construction, $\Hs_\ell\cap\Ext{S^3}{L}$
and $-\Hs_s\cap\Ext{S^3}{L}$ are both relative cobordisms
from $\page{f}{\pi/2}\cap\Ext{S^3}{L}$ to 
$\page{f}{-\pi/2}\cap\Ext{S^3}{L}$.
By the assumptions on $f(\crit f)$, 
$\page{f}{\pi/2}\cap\Ext{S^3}{L}$ and 
$\page{f}{-\pi/2}\cap\Ext{S^3}{L}$ are both 
diffeomorphic to $S_\ell$.  The claims follow immediately 
upon application of the usual dictionary between Morse functions 
and cobordisms (see \cite{Milnor:Morse}) to continuous branches
of $-\imunit\log(f\restr(\Hs_\ell\cap\Ext{S^3}{L}))$
and $-\imunit\log(f\restr(\Hs_s\cap\Ext{S^3}{L}))$.
\end{proof}

\begin{unproved}{cor}\label{from small to large--general}
$S_\ell$ is diffeomorphic to $S_s$ with $\nu$ hollow handles attached.
\end{unproved}

\begin{unproved}{cor}\label{from small to large}
The following are equivalent.  
\begin{inparaenum}
\item
$S_s$ is connected.
\item
$f$ is connected.
\item
$S_\ell$ is diffeomorphic to a connected sum
of $S_s$ and a closed surface of genus $\nu$.
\item
The interior of $\Ext{S^3}{L}$ contains arcs 
$\a_{+}, \a_{-}$ and handlebodies $\hb{\nu}_{+}, \hb{\nu}_{-}$
with the following properties.
\begin{inparaenum}
\item
$f(\hb{\nu}_{\pm})\sub\{\exp(\imunit\theta)\Suchthat 
|\theta\mp \pi/4|\le 1/10\}$.
\item
The arc $\a_{\pm}$ has one endpoint in $\Int S_s$ and the other in
$\Bd\hb{\nu}_{\pm}$, and is mapped diffeomorphically onto 
$\clint{0}{\pi/4 - 1/10}$ by the restriction of $\mp\imunit\log(f)$.
\item\label{ambient connected sum}
The ambient connected sum of $\Bd\hb{\nu}_{\pm}$ and $S_s$
along $\a_{\pm}$ is isotopic to $S_\ell$.  \textup(With an 
appropriate choice of the ``connecting tube'' which is
guided by $\a_{\pm}$, the isotopy can be achieved by
following the gradient flow of $\mp\imunit\log(f)$.\textup)
\end{inparaenum}
\end{inparaenum}
\end{unproved}
 
Let $S$ be a surface.  Say that a relative $1$--handle decomposition
of a $3$--manifold
\begin{equation*}\label{relative 1-handle decomposition}
M=(S\times\clint{-1}{1})\cup \bigcup_{k=1}^{n}\h1k,
\tag{$\ddagger$}
\end{equation*}
in which each $1$--handle $\h1k$ 
is attached either to $S\times \{-1\}$ or to $S\times \{1\}$,
is \bydef{involutized} by a map $\iota\from M\to M$ 
provided that $\iota$ is an orientation-reversing
involution such that 
$\iota(S\times\clint{-1}{1})=S\times\clint{-1}{1}$ and, for some $\e>0$,
$\iota\restr S\times \clint{-\e}{\e}\suchthat (x,t)\mapsto (x,-t)$.

\begin{unproved}{lem}\label{involutibility conditions}
\begin{inparaenum}
\item
If the decomposition \thetag{\ref{relative 1-handle decomposition}} 
can be involutized, then $n$ is even.
\item
If $n=0$ or if $S$ is connected, then the decomposition 
\thetag{\ref{relative 1-handle decomposition}} can be involutized.
\end{inparaenum}
\end{unproved}

When the relative $1$--handle decomposition \thetag{\ref{domain of eta-s}} 
of the domain of $\eta_s$, constructed (implicitly) in the course of the 
proof of Prop.~\ref{Hs and Hl as cobordisms} by applying Morse 
theory to $-\imunit\log(f)$, can be involutized, it may also
be said that $f$ can be involutized.  For example, if $f$ 
is connected (e.g., if $f$ is a fibration), then $f$ can be 
involutized; on the other hand, $u$ can be involutized 
although the small Seifert surface of $u$ is not connected.  
The Morse map $\omega$ described in 
Example~\ref{non-boundary-connected Morse map for U} 
cannot be involutized.

\begin{unproved}{cor}
If $f$ can be involutized, then 
$\Hs_\ell$ and $\Hs_s$ are handlebodies of genus $g$,
that is, $(\Hs_\ell,\Hs_s)$ is a genus-$g$ Heegaard 
splitting of $S^3$.
\end{unproved}

\begin{definitions}
\label{monodromy construction}  
\hypertarget{monodromy hyperconstruction}{}
Suppose $f$ can be involutized.
Choose diffeomorphisms $\eta_s$ and $\eta_\ell$, 
a relative handle decomposition \thetag{\ref{domain of eta-s}} 
of the domain of $\eta_s$, and an involution $\iota_s$ which 
involutizes $f$. The \bydef{monodromy} of $f$ 
(given these choices) is the diffeomorphism 
$h\isdefinedas \pr_1 \after \eta_\ell^{-1} \after \eta_s \after \iota_s %
\after \eta_s^{-1} \after \eta_\ell \after (\id{S_\ell},1) %
\from S_\ell\to S_\ell$ (pictured schematically 
in Fig.~\ref{schematic}). %
\begin{figure}
\centering
\includegraphics[width=.75\textwidth]{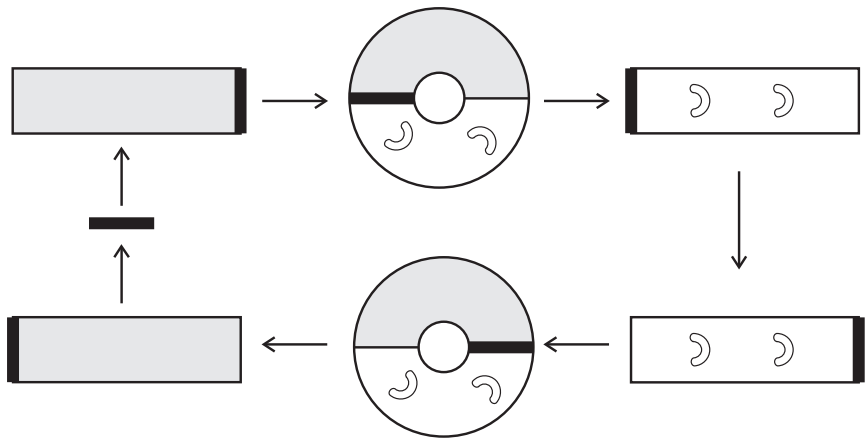}
\caption{Schematic decomposition of a monodromy.
\label{schematic}}
\end{figure}
Suppose further that $f$ is connected.
Choose diffeomorphisms $\d_{\pm}\from S_\ell\to S_s\connsum\hb{\nu}_{\pm}$
at the ends of isotopies as in 
Cor.~\ref{from small to large}\eqref{ambient connected sum}.
Choose a smooth map $c_{\pm}\from S_s\connsum\hb{\nu}_{\pm}\to S_s$ which 
(in an obvious way) ``collapses the second connected-summand to a point''.
The \bydef{$\pm$--adiexodon} of $f$ (given these choices) is the smooth 
map $A_{\pm}\isdefinedas c_{\pm}\after \d_{\pm}\from S_\ell\to S_s$.
\end{definitions}

\begin{figure}
\centering
\includegraphics[width=0.75\textwidth]{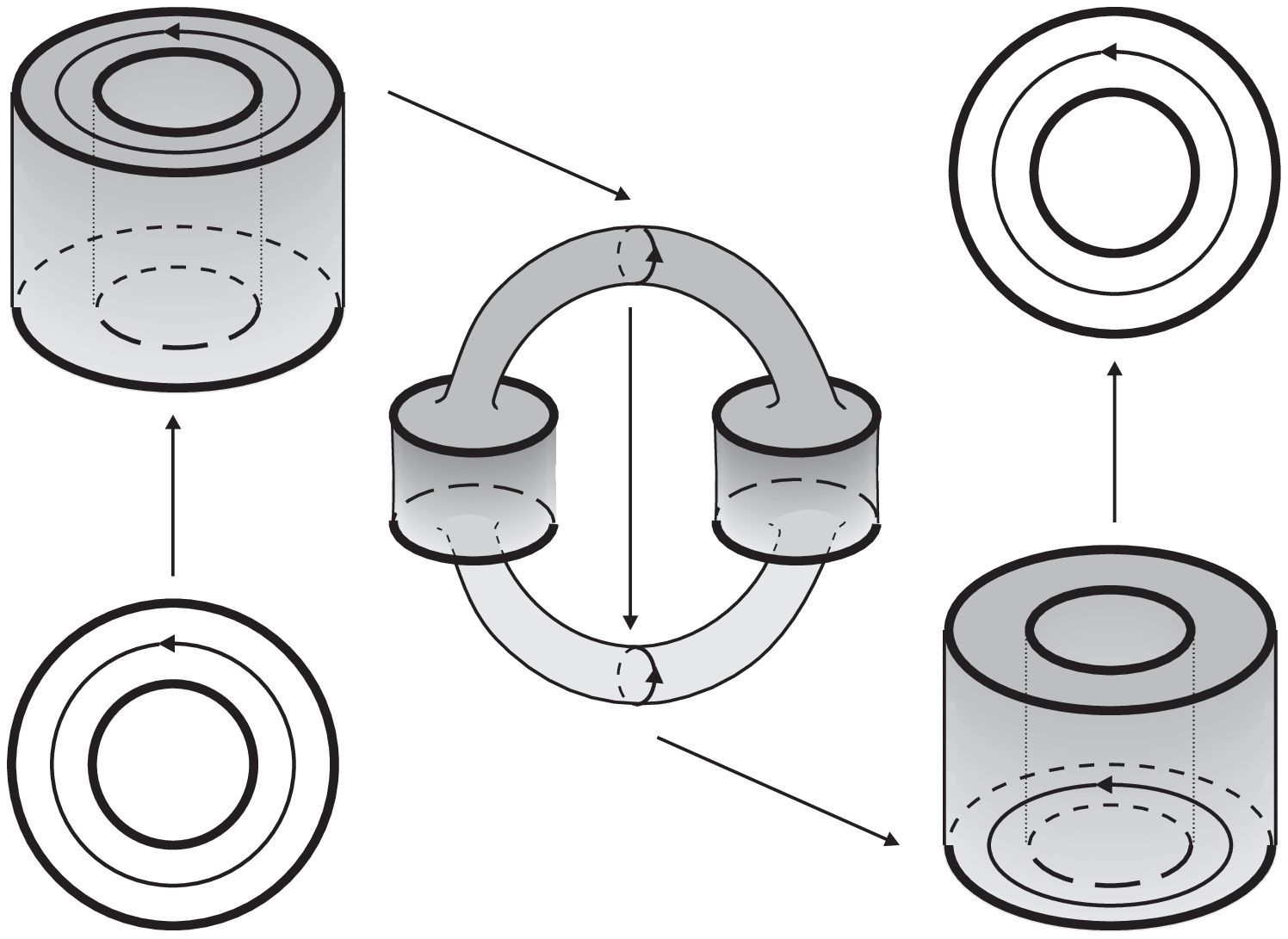}
\caption{\label{monodromy for u} A monodromy for $u$;
the middle arrow represents the composition 
$\eta_s\after\iota_s\after\eta_s^{-1}\from\Hs_s\cap\Ext{S^3}{U}
\to\Hs_s\cap\Ext{S^3}{U}$.}
\end{figure}

\begin{figure}
\centering
\includegraphics[width=0.9\textwidth]{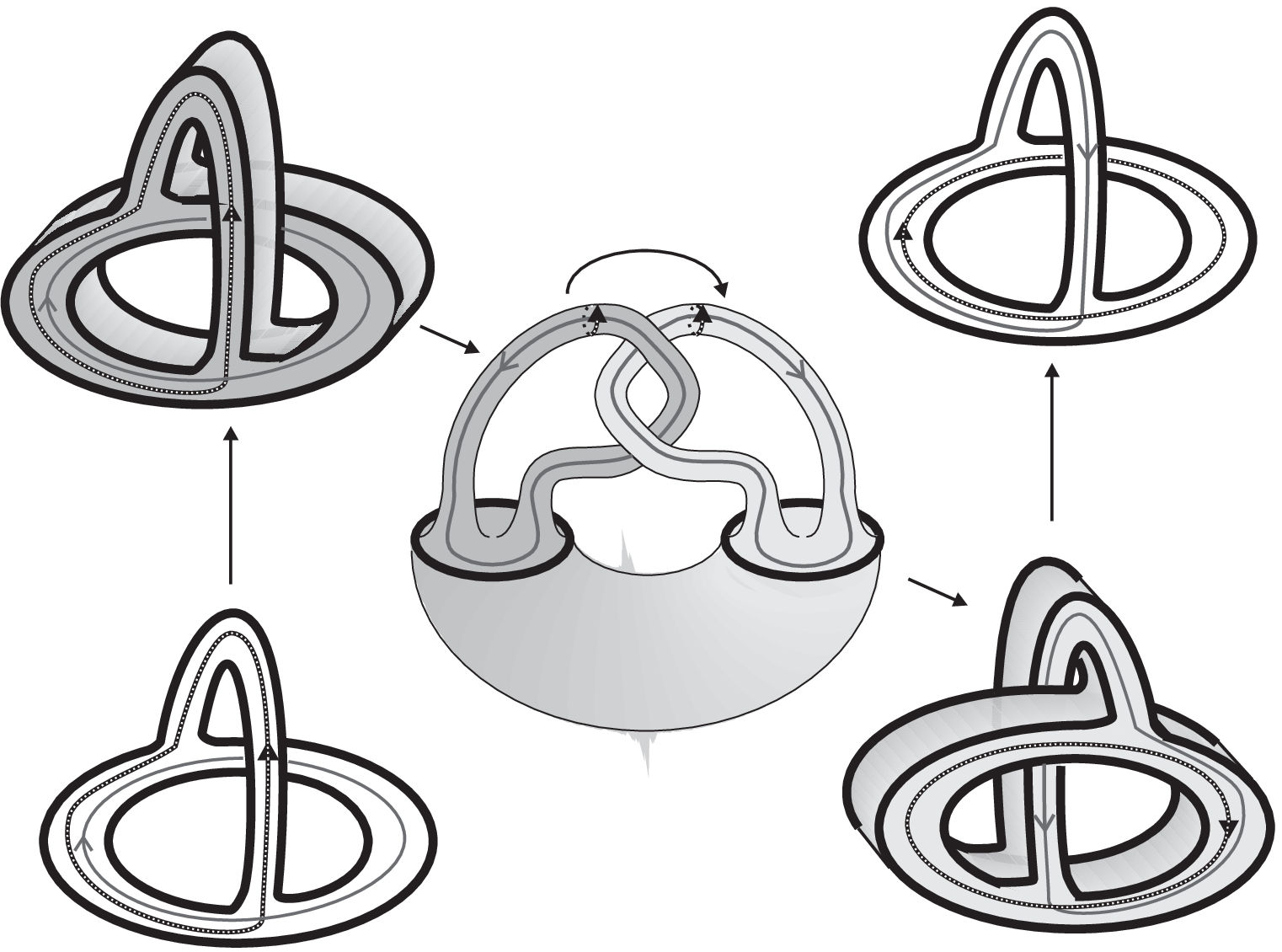}
\caption{\label{o1mono} A monodromy for $o_1$;
the middle arrow represents the composition 
$\eta_s\after\iota_s\after\eta_s^{-1}\from\Hs_s\cap\Ext{S^3}{O}
\to\Hs_s\cap\Ext{S^3}{O}$.}
\end{figure}

\begin{figure}
\centering
\includegraphics[width=0.9\textwidth]{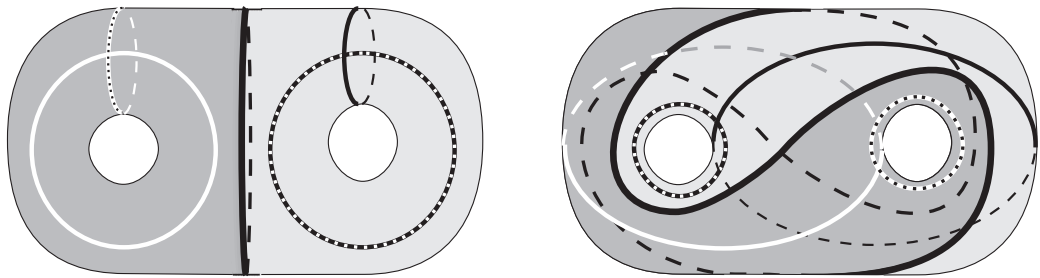}
\caption{\label{o1heeg} Abstract handlebodies 
diffeomorphic to the Heegaard handlebodies 
$\Hs_s, \Hs_\ell\sub S^3$ defined by 
$o_1\from\protect\diff{S^3}{O}\to S^1$.  On 
the surfaces bounding the handlebodies, 
corresponding loops are drawn in the same style;
in particular, the heavy, dark separating curves 
represent $O\sub S^3$.}
\end{figure}

A monodromy for $u$ (the identity map of an annulus)
is pictured in Fig.~\ref{monodromy for u}.
A monodromy for $o_1$ (in effect, the map 
$S^1\times S^1\suchthat (z,w)\mapsto(\conjug{w},z)$ 
restricted to $\Ext{S^1\times S^1}{\{(1,1)\}}$) 
is pictured in Fig.~\ref{o1mono}.
and a corresponding Heegaard splitting in Fig.~\ref{o1heeg}. %

\begin{remarks}
\begin{inparaenum}
\item
In case $f$ is a fibration, the monodromy of $f$ is a
geometric monodromy (or ``holonomy'', \cite{Stallings}) 
of $f$ in the usual sense, somewhat justifying so naming
it in general (but see note~\eqref{Katz} 
on page~\pageref{Katz});
both adiexodons of a fibration may be taken to be
the identity.  
\item\label{reconstruction}
It is well known that, when $f$ is a fibration,
both the embedding of $L$ into $S^3$ and the map $f$ 
can be reconstructed (up to isotopy) from $h$ alone, 
via the obvious construction of $\Ext{S^3}L$ 
as the mapping torus of $h$ (that is, the identification 
space $(S_\ell\times \clint{-1}{1}){/}{\equiv}_{h}$, 
where the non-trivial equivalence classes of ${\equiv}_{h}$ 
are the pairs $\{(x,-1),(h(x),1)\}$ for $x\in S_\ell$).  
In general, both the embedding of $L$ into $S^3$ and the map $f$ 
can be reconstructed (up to isotopy) similarly, using the extra 
data provided by the adiexodons.  In particular, all isotopy 
invariants of $L\sub S^3$ or $f$ are determined by 
the monodromy and adiexodons of $f$---though the 
calculation of any particular invariant may be more 
or less opaque (a familiar fact even when $f$ is a fibration).  
The Alexander invariant (in particular,
the Alexander polynomial) and Novikov homology, among others, 
do have reasonably transparent calculations in these terms.
Likewise, a suggestive realization of the kernel of 
$\pi_1(\diff{S^3}{L})\to \Z\suchthat [\g]\mapsto \int_\g \xi_L$
as an infinite free product with amalgamation, 
in the style of \cite{Neuwirth:algebraic} and
\cite{Brown-Crowell:augmentation}, 
falls out naturally from any monodromy and adiexodons of 
$f\from\diff{S^3}{L}\to S^1$.
However, we will not undertake such calculations in this paper.
\item
The triples $(\Hs_s,\page{f}{\pi/2},-\page{f}{-\pi/2})$ and
$(\Hs_\ell,-\page{f}{\pi/2},\page{f}{-\pi/2})$ are ``sutured manifolds''
as defined in \cite{Goda:suture-sum} (equivalent to
the original definition, \cite{Gabai:foliations-genera}).
The language of sutured manifolds could have been used
throughout this section, although the language of Morse maps to $S^1$
seems more natural for the study of monodromy and closed $f$-braids.
\end{inparaenum}
\end{remarks}

\begin{hists}\label{notes on monodromies}
\begin{inparaenum}
\item\label{Katz}
Of course the combination of Greek roots in ``monodromy''
means, more or less, ``single road'' (a denotation of 
\textsc{MONO$\Delta$POMO$\Sigma$} in modern Greek
is ``one-way street'').  As Nicholas Katz has pointed 
out to us \cite{Katz}, 
``In 19th century mathematics, monodromic meant
`single-valued'.  The `monodromy group'
of a differential equation measures the \emph{failure}
of its solutions to be monodromic, and by extension from 
that usage, `monodromy' has come to mean what should 
really be called polydromy.''
Since the word has, however unfortunately, become standard 
to describe the ``first return map'' from which a fibration 
over $S^1$ can be reconstructed, we feel somewhat 
justified in extending its use to the analogous part
of the geometric data from which a more general
Morse map from a link complement to $S^1$ can be 
reconstructed.
\item\label{Wyatt}
As noted after Definitions~\ref{monodromy construction}, to reconstruct
a Morse map which is not a fibration more is needed than
a monodromy.  In a sense, the extra information which is 
needed is carried by paths that, rather than ``going once 
around'' (whether to close up or keep going), go to---or come 
from---nowhere.  We are indebted to the classicist
William Wyatt \cite{Wyatt} for suggesting the word ``adiexodon'', 
coined by Aristotle in a discussion of the infinite \cite{Aristotle}, 
and parseable as (more or less) ``that from which there is no way 
out by going through''.
\end{inparaenum}
\end{hists}

\section{Splicing and closed $f$-braids}\label{splicing}

The theory of splicing (multi)links in homology $3$--spheres 
was introduced (for empty links) by Siebenmann \cite{Siebenmann:splice}, 
and extensively developed by Eisenbud \& Neumann \cite{Eisenbud-Neumann}.
We use only a special case, of which the following example is the paradigm.

\begin{example}\label{simplest splice example}
Let $L_0\isdefinedas L(G_0,\infty)\sub S^3$ 
be the fibered link of indeterminacy at infinity 
of the rational function $G_0\suchthat(z,w)\mapsto zw/(4z-1)$ 
introduced in Example~\ref{fibered loi@infinity}.  The three 
components of $L_0$ are unknots: the horizontal and vertical
unknots $O$ and $O'$ (as on p.~\pageref{unknots definition}),
and $O_0\isdefinedas\{(z,w)\in S^3\Suchthat z=1/4\}$. 
With the orientation imposed on $L_0$ by $\Mmap{G_0}$, 
$O'-O_0$ is the oriented boundary of an annulus 
$\AKn{O'}{0}\sub\diff{S^3}{O}$.  This fact 
(or direct calculation) shows that, if 
$\t_0\from L_0\times D^2 \to \Nb{S^3}{L_0}$ is a 
trivialization which is adapted to $\Mmap{G_0}$, then 
$\t_0\restr O\times D^2$ is isotopic (by an 
isotopy that can be made arbitrarily small by
shrinking the regular neighbhorhoods) to
a trivialization $\t\from O\times D^2\to\Nb{S^3}{O}$
which is adapted to $o$.  Of course the exterior 
$\Ext{S^3}{O}$ of $O$ is also a regular neighborhood 
$\Nb{S^3}{O'}$ of $O'$.  From that point of view, 
up to an arbitrarily small isotopy $L_0$ and $\Mmap{G_0}$ 
can be thought of as being constructed from $O$ 
and $o$ by replacing the triple 
\[
(\Nb{S^3}{O'},O',o\restr(\diff{\Nb{S^3}{O'}}{O'}))
\]
with the triple
\[
(\Ext{S^3}{O},O'\cup O_0,\Mmap{G_0}\restr(\diff{\Ext{S^3}{O}}{O'\cup O_0})).
\]
At the level of a single fiber surface, this can be done
in two steps:
first, replace a $2$-disk $\page{o}{\theta}$
by $\page{o}{\theta}\cup\AKn{O'}{0}$, a \bydef{Seifert ribbon} 
(see \cite{braided-surfaces}) bounded by $L_0$, immersed with 
doublepoints along the arc $\page{o}{\theta}\cap\AKn{O'}{0}$); 
then disingularize $\page{o}{\theta}\cup\AKn{O'}{0}$
by ``opening up'' its arc of doublepoints into an unknot.
(The situation is illustrated in Fig.~\ref{simplest splice}.)
\begin{figure}
\centering
\includegraphics[width=0.7\textwidth]{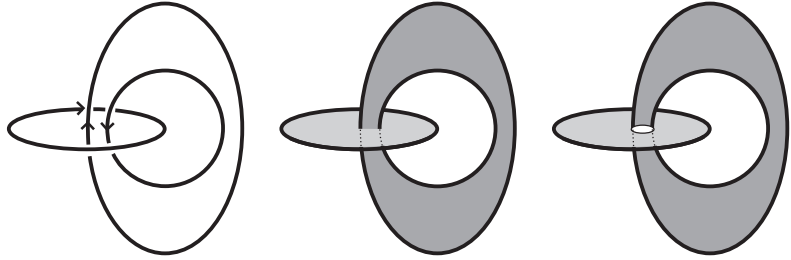}
\caption{The link $L_0$; a Seifert ribbon for $L_0$;
a fiber surface for $L_0$.
\label{simplest splice}}
\end{figure}
Actually, as the construction shows, this can be done for all
fiber surfaces simultaneously---a striking, though elementary, fact.

Entirely similarly, for all $k$ the fibered link of 
indeterminacy at infinity $L_k\isdefinedas L(G_k,\infty)\sub S^3$ 
of the rational function $G_k\suchthat(z,w)\mapsto zw/(4z-w^k)$ 
(also introduced in Example~\ref{fibered loi@infinity}),
and its fibration $\Mmap{G_k}$, can be constructed 
from $O$ and $o$ by replacing
\[
(\Nb{S^3}{O'},O',o\restr(\diff{\Nb{S^3}{O'}}{O'}))
\]
with
\[
(\Ext{S^3}{O},O'\cup O_k,\Mmap{G_k}\restr(\diff{\Ext{S^3}{O}}{O'\cup O_k})),
\]
where now $O_k\isdefinedas \{(z,w)\in S^3\Suchthat z=w^k/4\}=
\diff{L_k}{(O\cup O')}$, so that $O'-O_k$ is the oriented 
boundary of an annulus $\AKn{O'}{k}\sub\diff{S^3}{O}$; again, 
each fiber surface is constructed from a Seifert ribbon by
opening up an arc of doublepoints into an unknot.
\end{example}

What makes this example work is that $O'$ is a
closed braid around the axis $O$.  More precisely,
$O'$ is a closed $1$-string $o$-braid in the sense
of the following definition.

\begin{definition} 
Let $f\from \diff{S^3}{L} \to S^1$ be a Morse map.
If $B\sub\diff{S^3}{L}$ is a link such that
\begin{inparaenum}[(a)]
\item\label{f-braids miss crit(f)}
$B\cap\crit{f}=\emptyset$,
\item\label{f-braids are transverse}
$B$ intersects every page of $f$ transversely, and
\item\label{f-braids are positive}
every intersecton of $B$ with a page of $f$ is positive,
\end{inparaenum}
then $f\restr B\from B\to S^1$ is a local diffeomorphism
and an orientation-preserving covering projection of some 
degree $n\ge 1$.  Call such a link $B$ a 
\bydef{closed $n$-string $f$-braid}.  
Call $B$ \bydef{pure} if it has the same number 
of components as it has strings, that is, if each
component of $B$ is a closed $1$-string $f$-braid.  
\end{definition}

\begin{construction}[splicing along a closed $f$-braid] 
\label{splicing construction}
\hypertarget{splicing hyperconstruction}{}
Let $L\sub S^3$ be a link, 
$f\from\diff{S^3}{L}\to S^1$ a Morse map.
A \bydef{framed closed $n$-string $f$-braid} 
is a pair $(B,\framing{k})$ where $B$ is a 
closed $n$-string $f$-braid and $\framing{k}\from B\to\Z$ 
is a \bydef{framing} of $B$ (i.e., $\framing{k}$
is continuous).  As a closed $f$-braid, $B$ has 
a regular neighborhood with a trivialization 
$\t\from B\times D^2\to\Nb{\diff{S^3}{L}}{B}$
such that
\begin{inparaenum}[(a)]
\item\label{splicetriv1}
$\t(x,0)=x$ and $f(\t(x,z))=f(x)$ 
for all $(x,z)\in B\times D^2$, and
\item\label{splicetriv2}
for each component $C\sub B$, 
$\t(C\times\clint{0}{1/2})$
is an annulus $\AKn{C}{\framing{k}(C)}$.
\end{inparaenum}
Given such a trivialization $\t$, orient $\t(C\times\clint{0}{1/2})$
so that its boundary contains $C$ (rather than $-C$),
and let the \bydef{spliced link of $(B,\framing{k})$ along $f$}
be
\[
\splicedlink{f}{B}{\framing{k}}\isdefinedas 
L\cup\bigcup_{C\sub B} \Bd \t(C\times\clint{0}{1/2})
\]
Clearly, $\splicedlink{f}{B}{\framing{k}}$ depends (up to isotopy)
only on the embeddings of $L$ in $S^3$ and $B$ in $\diff{S^3}{L}$
(up to isotopy) and the framing $\framing{k}$, and not on $f$ 
(except insofar as $B$ is restricted to be a closed $f$-braid) 
nor on the choices of $\Nb{\diff{S^3}{L}}{B}$ and $\t$.  
Equally clearly, 
a Seifert surface for $\splicedlink{f}{B}{\framing{k}}$ can
be constructed from any smooth page $\page{f}{\theta}$, 
by taking $\page{f}{\theta}\cup \t(C\times\clint{0}{1/2})$
(a Seifert ribbon for $\splicedlink{f}{B}{\framing{k}}$ 
with $n$ arcs of doublepoints, the components of 
$\page{f}{\theta}\cap\t(C\times\clint{0}{1/2})$) and
open up each arc of doublepoints into an unknot 
as in Example~\ref{simplest splice example}.
Up to isotopy, this Seifert surface depends only on 
$\page{f}{\theta}$, $B$, and $\framing{k}$.

The \bydef{spliced map of $(B,\framing{k})$ along $f$}, 
$\splicedmap{f}{B}{\framing{k}}\from%
\diff{S^3}{\splicedlink{f}{B}{\framing{k}}}\to S^1$,
depends (up to isotopy) on $f$ as well as on $L$, $B$, 
and $\framing{k}$.  In case $B$ is pure, we mimic
Example~\ref{simplest splice example} exactly for
each component $C$ of $B$, replacing
\[
(\Nb{S^3}{C},C,f\restr(\diff{\Nb{S^3}{C}}{C}))
\]
with
\[
(\Ext{S^3}{O},O'\cup O_k,\Mmap{G_0}\restr(\diff{\Ext{S^3}{O}}{O'\cup O_k})).
\]
For general $B$, if the component $C\sub B$ is a closed
$\nu$-string $f$-braid, then 
\[
\Nb{S^3}{C},C,f\restr(\diff{\Nb{S^3}{C}}{C}))
\]
is replaced with
\[
(\Ext{S^3}{O},O'\cup O_{\framing{k}(C)},%
\Mmap{G_0}^\nu\restr(\diff{\Ext{S^3}{O}}{O'\cup O_{\framing{k}(C)}})).
\]
For more details (from a rather different, much more general,
viewpoint), see \cite{Eisenbud-Neumann}.
\end{construction}

\begin{example}\label{changing framing by twisting}
For all $f$ and $B$, and any framings $\framing{k}$, 
$\framing{k'}$ of $B$, $\splicedmap{f}{B}{\framing{k'}}$ 
is a composite of Stallings twists of 
$\splicedmap{f}{B}{\framing{k}}$ along the unknots
produced by opening up the arcs of doublepoints on
$\page{f}{\theta}\cup \t(C\times\clint{0}{1/2})$.
(In light of Example~\ref{simplest splice example},
this is a generalization of Example~\ref{fibered loi@infinity, bis}.)
\end{example}

\begin{unproved}{prop}\label{critical points of spliced map} 
For all $f$, $B$, and $\framing{k}$, 
$\crit{\splicedmap{f}{B}{\framing{k}}}=\crit{f}$,
and $\ind{\splicedmap{f}{B}{\framing{k}}}{x}%
\linebreak[1]=\linebreak[1]\ind{f}{x}$ 
for all $x\in\crit f$.  In particular, if $f$ is a fibration 
\textup(resp., moderate; self-indexed\textup), then so is 
$\splicedmap{f}{B}{\framing{k}}$.
\end{unproved}

\begin{unproved}{cor}\label{splicing doesn't make MN bigger} 
For all $L$, 
for every Morse map $f\from\diff{S^3}{L}\to S^1$ and 
every framed $f$-braid $(B,\framing{k})$, 
$\MN(\splicedmap{f}{B}{\framing{k}})\le \MN(L)$.
\end{unproved}

To make good use of Construction~\ref{splicing construction}, 
we need an extensive supply of closed $f$-braids.  For the 
rest of this section, assume that $f$ is moderate, self-indexed, 
and connected.  In this case, closed $n$-string $f$-braids are 
especially easy to describe---and, somewhat surprisingly, plentiful
and diverse even for $n=1$ unless $f$ is isotopic to $o$.

Let $f\from \diff{S^3}{L} \to S^1$ be normalized 
as at the beginning of Sect.~\ref{monodromies}.
Fix $n\ge1$ and $X\sub\Int\Nb{S_s}{\Bd S_s}$ with $\card{X}=n$.

\begin{prop}\label{description of closed f-braids}
If $B$ is a closed $n$-string $f$-braid, then, up to isotopy
through closed $n$-string $f$-braids:
\begin{enumerate}
\thmlistfix
\item\label{braid exterior to binding}
$B\sub\Ext{S^3}{L}$;
\item\label{braid in Hl}
$B\cap\Hs_\ell=\eta_\ell(B_\ell)$, where 
$\pr_2\restr B_\ell\from B_\ell\to\clint{-1}{1}$ 
is an orientation-prserving covering projection 
\textup(i.e.,
$B_\ell\sub S_\ell\times \clint{-1}{1}$ is a
\textup{geometric $n$-string $S_\ell$-braid}\textup); %
\item\label{braid in Hs}
\begin{inparaenum}
\item\label{braid misses 1-handles}
$B\cap\Hs_s\sub \eta_s(S_s\times \clint{-1}{1})$, and in fact
\item\label{braid is a product}
$B\cap\Hs_s=\eta_s(X\times \clint{-1}{1})$.
\end{inparaenum}
\end{enumerate}
Conversely, if $B$ has properties 
\eqref{braid exterior to binding}, \eqref{braid in Hl},
and \eqref{braid in Hs}, then $B$ is a closed $n$-string $f$-braid.
\end{prop}

\begin{proof} 
Suppose $B$ is a closed $n$-string $f$-braid.
Since $B\sub\diff{S^3}{L}$ is compact, certainly 
$B\cap N=\emptyset$ for some regular neighborhood 
$N$ of $L$ in $S^3$ which is concentric to $\Nb{S^3}{L}$
with respect to some trivialization of $\Nb{S^3}{L}$ that
is adapted to $f$.  An appropriate isotopy of 
$S^3$ (rel.\ $L$) supported near $\Nb{S^3}{L}$, 
preserving level sets of $f$ and thus
the class of closed $n$-string $f$-braids, 
carries $N$ onto $\Nb{S^3}{L}$ and $B$ onto 
a closed $n$-string $f$-braid with 
property~\eqref{braid exterior to binding}.
Now property~\eqref{braid in Hl} is immediate 
from Prop.~\ref{Hs and Hl as cobordisms}\eqref{f on Hl}.

\begin{figure}
\centering
\includegraphics[width=\textwidth]{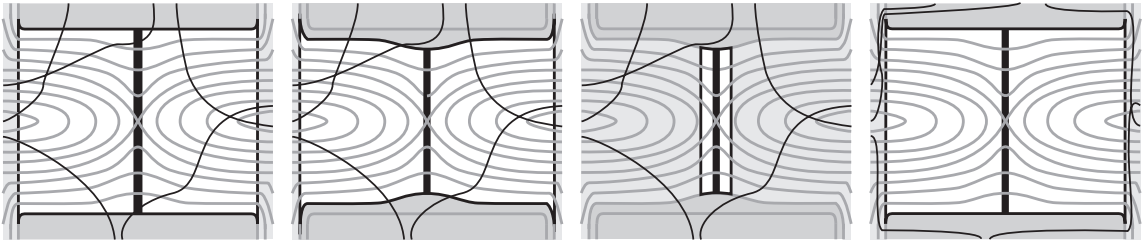}
\caption{Each frame shows 
part of \protect$\Hs_\ell$ (dark gray),
\protect$\eta_s(S_s\times\clint{-1}{1})$ (light gray),
and $B$ (black curves), level sets of $f$ (gray curves), 
and all of \protect$\eta_s(\h1{k+r\nu})$ (white with black boundary)
and \protect$\D_{k+r\nu}$ (thick black line).
\label{handle figure}}
\end{figure}
Suppose $B$ has property~\eqref{braid exterior to binding}.
From the construction of $\eta_s$ and the handle 
decomposition \thetag{\ref{domain of eta-s}} of its domain, 
for $k=1,\dots,\nu$ and $r=0,1$ there exists a transverse 
$2$--disk $\D_{k+r\nu}\sub \eta_s(\h1{k+r\nu})\sub\Hs_s\cap\Ext{S^3}{L}$
of the handlebody $\Hs_s$ such that $f\restr \D_{k+r\nu}$ is 
Morse with exactly one critical point $x_{k+r\nu}$, of index $2r$,
and $\eta_s(\h1{k+r\nu})\cap\crit{f}=%
\D_{k+r\nu}\cap f^{-1}(\exp((-1)^{r+1}\imunit\pi/4)=%
\{\x_{k+r\nu}\}$.   By \eqref{f-braids miss crit(f)} 
and \eqref{f-braids are transverse},
the finite set $B\cap \bigcup_{k,r}\D_{k+r\nu}$ is disjoint
from $\crit{f\restr\bigcup_{k,r}\D_{k+r\nu}}$.
Now the process adequately sketched in Fig.~\ref{handle figure}
(swapping a collar of their common boundary from one handlebody 
to the other to shrink the transverse disks and make
$B\cap \bigcup_{k,r}\D_{k+r\nu}=\emptyset$; 
shortening the $1$--handles around their shrunken transverse disks 
to make $B\cap \bigcup_{k,r}\eta_s(\h1{k+r\nu})=\emptyset$; 
isotoping the handlebody decomposition with shrunken,
shortened $1$--handles back to the original one, 
and carrying $B$ along by the isotopy) lets
$B$ acquire property~\eqref{braid misses 1-handles} while
preserving properties \eqref{braid exterior to binding} and
\eqref{braid in Hl}.
Since $S_s$ is connected, any two embeddings of the same
finite set into $\Int{S_s}$ are isotopic, so it is easy to
endow $B$ with the further property~\eqref{braid is a product}
while preserving properties 
\eqref{braid exterior to binding}, \eqref{braid in Hl},
and \eqref{braid misses 1-handles}.

The converse is trivial.
\end{proof}

\begin{cor}\label{knot on page gives 1-string braid}
If a knot $K$ is contained in a page of $f$, 
then $K$ is isotopic in $S^3$ to a closed $1$-string 
$f$-braid $B$.
\end{cor}

\begin{proof}
Let a link $K$ be contained in the page $\page{f}{\theta}$ 
of $f$. Let $x\in K$ be a basepoint.

If $\page{f}{\theta}$ is singular and $r$ is the index of 
$f$ at any (and therefore every) point of $\crit{f}\cap \page{f}{\theta}$, 
then (because $r\in\{1,2\}$) it is easy to construct a 
vectorfield $V$ on a small neighborhood of $K$ such that 
$V$ is non-zero off $L$, $V$ equals $(-1)^{r+1}$ times 
the gradient of $-\imunit\log(f)$ outside a smaller 
neighborhood of $\crit{f}\cap \page{f}{\theta}$, and the 
flow of $V$ carries $K$ onto a smooth page of $f$.  

Suppose $\page{f}{\theta}$ is a smooth page of $f$.
There is an arbitrarily small isotopy of $S^3$ 
which carries $\page{f}{\theta}$ (and therefore $K$) 
into $\Int \page{f}{\theta}$.  By another (possibly lengthy)
isotopy, $K$ can be moved into $S_\ell$ and $x$ into 
$\Int\Nb{S_\ell}{\Bd S_\ell}$.

Now suppose that $K$ is a knot.  
Let $x_s\in\Int\Nb{S_\ell}{\Bd S_\ell}$ be
the point such that $\eta_s(x_s,\pm1)=\eta_\ell(x,\mp1)$.
Let $\g\from(\clint{-1}{1},1)\to(K,x)$
be a $\mathscr C^\infty$ map onto $K$ such that 
$\g\restr\opint{-1}{1}$ is a bijective immersion
and $\g$ is infinitely flat at $\pm 1$.  
If $B_\ell\sub S_\ell\times\clint{-1}{1}$ 
is the geometric $1$-string $S_\ell$ braid parametrized by 
$t\mapsto (\g(t),t)$, then 
$B\isdefinedas \eta_\ell(B_\ell)\cup \eta_s(\{x_s\}\times\clint{-1}{1})$
is a closed $1$-string $f$-braid isotopic to $K$ in $S^3$. %
\begin{figure}
\centering
\includegraphics[width=0.75\textwidth]{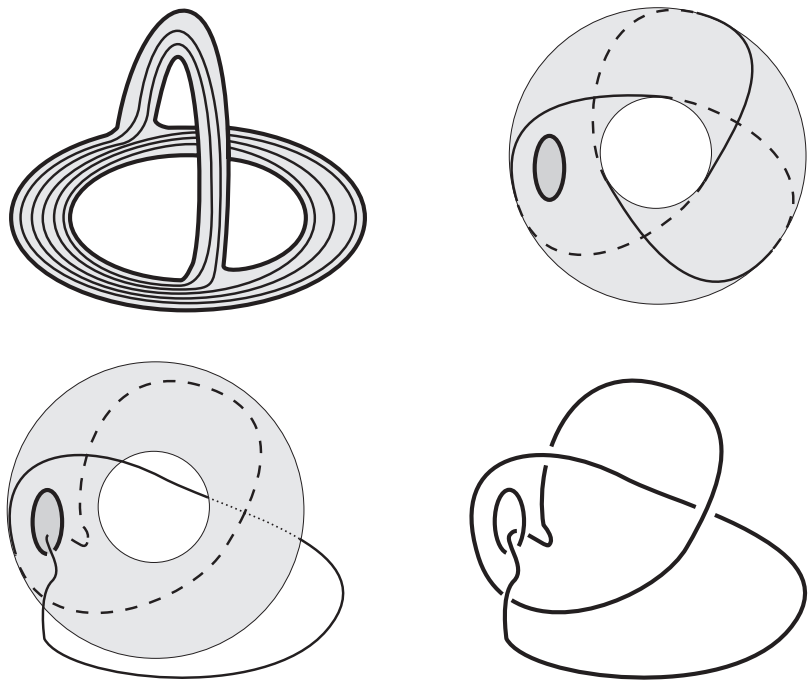}
\caption{Top: $K$, a torus knot of type $(2,3)$ 
lying on the large Seifert surface of $o_1$
(left, \protect$\page{o_1}{-\pi/2}$ as pictured in Fig.~\protect\ref{o1mono}; 
right, \protect$\page{o_1}{-\pi/2}$ with its unknotted boundary in standard
position, together with the disk $\protect\page{o_1}{\pi/2}$).
Bottom: a closed $1$-string $o_1$-braid $B$ isotopic to $K$
in $S^3$ (left, the part of $B$ inside the torus
\protect$\page{o_1}{\pi/2}\cup \page{o_1}{\pi}$ 
is represented by dashed lines,
and the rest of $B$ by solid and dotted lines; 
right, all of the torus but $O$ has been removed).
\label{O{2,3} as a 1-string o1-braid}}
\end{figure}
(The situation is illustrated in Fig.~\ref{O{2,3} as a 1-string o1-braid}.)
\end{proof}

\begin{figure}
\centering
\includegraphics[width=0.75\textwidth]{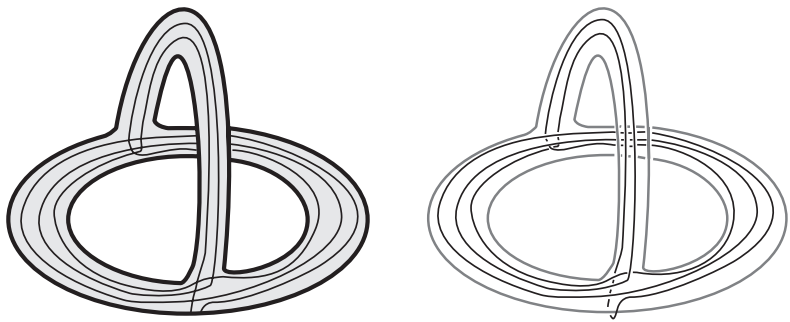}
\caption{Left: an $o_1$-ascending knot diagram with
eight crossings, on a large Seifert surface of $o_1$.  
Right: the union of the $1$-string $o_1$-braid with 
the given diagram, and the binding $O$ of $o_1$; this
link is isotopic to the union of a figure-$8$ knot and 
one of its meridional circles.
\label{figure-8 layered on torus}} 
\end{figure}

Prop.~\ref{description of closed f-braids}\eqref{braid in Hl}
shows immediately that (up to isotopy) any closed $1$-string 
$f$-braid may be constructed in a similar manner, starting 
from an immersed (not necessarily embedded) proper arc $\g$
on a smooth page $\page{f}{\theta}$ and treating $\g$ as an 
``$f$--ascending knot diagram'' on $\page{f}{\theta}$.  
This construction is illustrated in 
Fig.~\ref{figure-8 layered on torus}.

The \bydef{wrapping genus} $\wrapgenus{K}$
of a knot $K$ has been defined as the least $n$ 
such that $K$ embeds on a Heegaard surface of 
genus $n$ (so, for instance, $\wrapgenus{K}=1$
if and only if $K$ is a non-trivial torus knot).  
Since a large Seifert surface of the $g$-fold 
connected sum $o_g\isdefinedas \Connsum_1^g o_1$ 
is the exterior of a point on a Heegaard surface 
of genus $g$, it is equivalent to define 
$\wrapgenus{K}$ as the least $g$ such that, 
up to isotopy, $K$ is embedded on a large 
Seifert surface of $o_g$.  Let $\wlapgenus{K}$,
the \bydef{layered wrapping genus} of $K$,
be the least $g$ such that $K$ is isotopic 
to a closed $1$-string $o_g$-braid.  (The name
is appropriate in view of the description of 
closed $1$-string $f$-braids in the previous paragraph.)

\begin{unproved}{cor}\label{wlap at most wrap}
For all $K$, $\wlapgenus{K}\le\wrapgenus{K}$.
\end{unproved}

Since the figure-$8$ knot is not a torus knot,
Fig.~\ref{figure-8 layered on torus} shows that
the inequality in Cor.~\ref{wlap at most wrap}
can be sharp.

\begin{cor}\label{first wrapping estimate}
If $K$ is a knot and $D\sub S^3$ is a $2$--disk 
such that $D\cap\AKn{K}{k}$ is a transverse arc
of an annulus $\AKn{K}{k}$, then 
$\MN(\Bd D\cup\Bd\AKn{K}{k})\le 2\wlapgenus{K}$
for all $k\in\Z$.
\end{cor}
\begin{proof} 
We may assume that $K$ is a closed 
$1$-string $o_{\wlapgenus{K}}$-braid.
By Cor.~\ref{knot on page gives 1-string braid},
Construction~\ref{splicing construction},
and Prop.~\ref{critical points of spliced map},
$\card{\crit{\splicedmap{o_{\wlapgenus{K}}}{K}{k}}}=
\card{\crit{o_{\wlapgenus{K}}}}=2\wlapgenus{K}$.
But $\splicedlink{o_{\wlapgenus{K}}}{K}{k}$ and 
$\Bd D\cup\Bd\AKn{K}{k}$ are isotopic.
\end{proof}

\begin{remarks}
\begin{inparaenum}
\item
We thank Claude Weber for pointing out to us that 
Construction~\ref{splicing construction} is an instance 
of splicing.
\item
The first author (re)discovered Construction~\ref{splicing construction},
for closed $o$\,-braids (using geometric techniques quite different
than those described in this paper), as part of the authors' on-going joint 
investigations into Hopf plumbing and the geometry of the enhanced 
Milnor number, \cite{Hirasawa-Rudolph:enhanced}; cf.~\cite{icp1}, 
\cite{Giroux:ICM}.
\item
Eisenbud \& Neumann \cite[Theorem 4.2]{Eisenbud-Neumann} 
prove that a (multi)link is fibered if and only if it is 
non-split (``irreducible'') and each of its (multi)link 
``splice components'' is fibered.  It is clear that the 
methods of \cite{Eisenbud-Neumann} could be used to
derive quite general (sub)additivity formul{\ae} for 
Morse--Novikov number over appropriate splicings
(including cablings).  
\end{inparaenum}
\end{remarks}

\begin{hist} 
\begin{inparaenum}
\item
A closed $o$\,-braid is simply a closed braid
with axis $O$ in the usual sense.  There is 
some published work \cite{constqp1,Skora,Sundheim}
on closed $f$-braids for a fibration $f\from \diff{S^3}{L}\to S^1$ 
(and, more generally, for a fibration $f\from M\to S^1$ 
where the $3$--manifold $M$ is not necessarily a link complement 
in $S^3$).  
\item
Stillwell \cite{Stillwell} calls $\wrapgenus{K}$ the
``handle number'' of $K$; Goda's ``handle number'' 
\cite{Goda:handle-number} is something else entirely
(see the note on p.~\pageref{Goda's priority}).
\item
We are not aware of previous work on closed 
$f$-braids for Morse maps $f$ which are not fibrations,
nor on layered wrapping genus.
\end{inparaenum}
\end{hist}

\section{Murasugi sums of Morse maps}\label{Murasugi sums}

The underlying idea of the construction of a 
$2n$--gonal Murasugi sum $f_0\plumb{(N_0,N_1,\zeta)}f_1$ 
of two Morse maps $f_s\from \diff{S^3_s}{L_s}\to S^1$ 
is simple: given appropriate $2n$--gonal $2$--cells 
(so-called $n$--patches) on a smooth page of $f_0$ 
and a smooth page of $f_1$, with appropriate regular
neighborhoods $N_s\sub S^3_s$ (each a $3$--disk) 
such that $f_s\restr N_s$ is ``standard'' in a suitable 
sense, and in particular has no critical points in $\Int N_s$, 
the exteriors $E_s=\diff{S^3}{\Int(N_s)}$ (which are also
$3$--disks) may be glued together along their boundaries
to produce a $3$--sphere on which the restrictions
$f_s\restr E_s$ glue together to give a Morse map
(on the complement of an appropriate link, denoted
$f_0\plumb{(N_0,N_1,\zeta)}f_1$) whose critical points 
are the disjoint union of the critical points of $f_0$ 
and $f_1$.  The details, unfortunately, are somewhat 
technical.

For $n\in\N$, write $\G_{n}\isdefinedas \{z\in\C\Suchthat z^{n}=1\}$.
Let $p_n\from \diff{\Cext}{\G_{2n}}\to S^1$
be the argument of the rational function
$\Cext\to\Cext\suchthat z\mapsto(1+z^n)/(1-z^n)$,
viz., 
\[
p_n(z)=\frac{(1+z^n)/|1+z^n|}{(1-z^n)/|1-z^n|} 
\textrm{ for $z\notin \G_{2n}\cup\{\infty\}$, }
p_n(\infty)=-1.
\]
Define $P_n\from \diff{\Cext}{\G_{2n}}\times\clint{0}{\pi}\to S^1$ 
by $P_n(z,\theta)=\exp(\imunit\theta)p_n(z)$.  Note that
$p_n(\zeta z)=p_n(z)$ and $P_n(z,\theta)=P_n(\zeta z,\theta)$
for any $\zeta\in \G_n$.

\begin{unproved}{lem}\label{description of pn and Pn}
\begin{inparaenum}
\item
$p_1$ is a fibration with fiber 
$\opint{-1}{1}=p_1^{-1}(1)$.
For $n>1$, $p_n$ has exactly two critical points, 
$0\in p_n^{-1}(1)$ and $\infty\in p_n^{-1}(-1)$,
at each of which the germ of $p_n$ is smoothly 
conjugate to the germ of $z\mapsto\Im(z^n)$ at $0$.
\item
$P_1$ is a fibration with fiber 
$\opint{-1}{1}\times\clint{0}{\pi} =P_1^{-1}(1)$.
For all $n$, $P_n$ has no critical points,
and there is a trivialization 
$\t_{P_n}\from (\G_{2n}\times\clint{0}{\pi})\times D^2%
\to\Nb{\Cext\times\clint{0}{\pi}}{\G_{2n}\times\clint{0}{\pi}}$ 
which is adapted to $P_n$.
\end{inparaenum}
\end{unproved}

\begin{figure}
\centering
\includegraphics[width=0.6\textwidth]{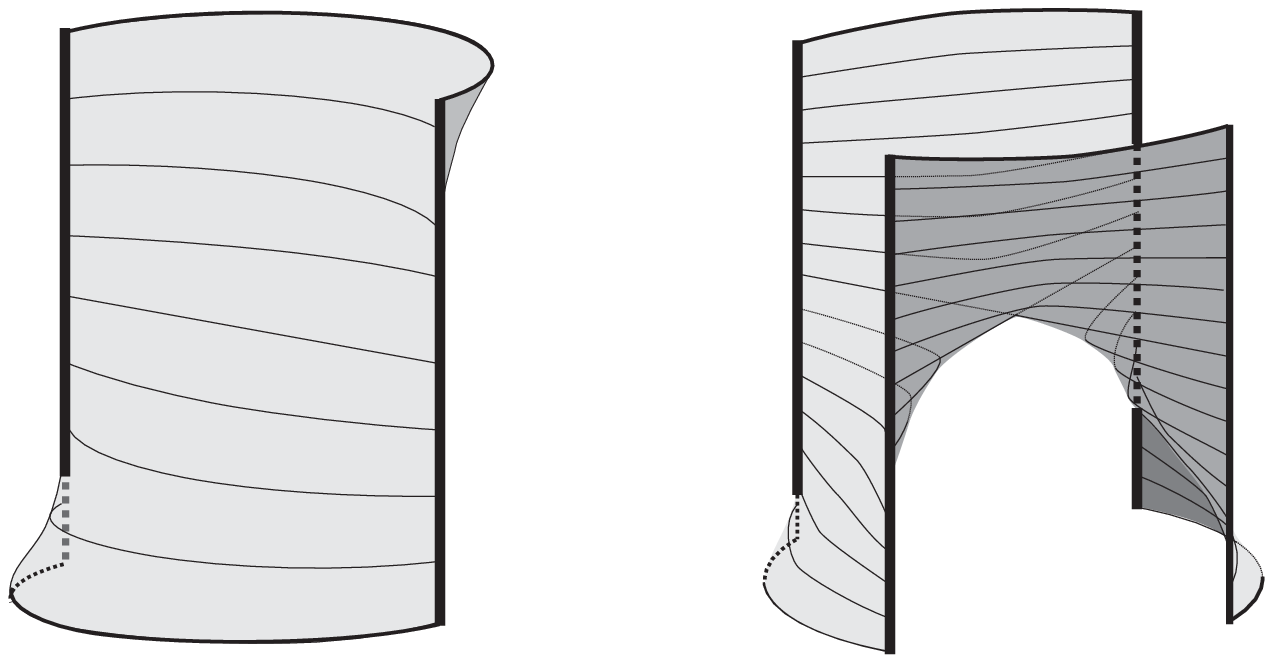}
\caption{$Q_1(\pi/2)$ and some 
level sets of $\pr_2\protect\restr Q_1(\pi/2)$; 
$Q_2(\pi/2)$ and some level sets of $\pr_2\protect\restr Q_2(\pi/2)$.
\label{q1&q2}}
\end{figure}

Let $Q_n(\theta)$ denote the closure of $P_n^{-1}(\exp(\imunit\theta))$ 
in $\Cext\times\clint{0}{\pi}$.  
By inspection (and Lemma~\ref{description of pn and Pn}), 
for all $\exp(\imunit\theta)\in S^1$, $Q_1(\theta)$ is a $4$--gonal $2$--disk
with smooth interior bounded by the union of 
$\G_{2}\times\clint{0}{\pi}$ and $2$ semicircles in $\Cext\times\{0,\pi\}$;
for $n>1$ and all $\exp(\imunit\theta)\in \diff{S^1}{\{1,-1\}}$, 
$Q_n(\theta)$ is a $4n$--gonal $2$--disk
with smooth interior bounded by the union of
$\G_{2n}\times\clint{0}{\pi}$ and $2n$ circular arcs 
in $\Cext\times\{0,\pi\}$.  
(The cases $n=1, 2$, $\theta=\pi/2$, are pictured in 
Fig.~\ref{q1&q2}.)

\begin{unproved}{lem}
\label{description of Qn}
For all $\exp(\imunit\theta)\in S^1$, the restriction 
$\pr_2\restr Q_1(\theta)\from Q_1(\theta)\to\clint{0}{\pi}$ has 
no critical points.  
For $n>1$ and $\exp(\imunit\theta)\in \diff{S^1}{\{1,-1\}}$, 
there is exactly one critical point 
of $\pr_2\restr Q_n(\theta)\from Q_n(\theta)\to\clint{0}{\pi}$ 
\textup{(}to wit,
$(0,\theta)$ for $0<\theta<\pi$
and $(\infty,\theta-\pi)$ for $\pi<\theta<2\pi$\textup{)},
at which the germ 
of $\pr_2\restr Q_n(\theta)$ is smoothly conjugate to the germ of 
$z\mapsto\Im(z^n)$ at $0$.
\end{unproved}

\begin{figure}
\centering
\includegraphics[width=0.6\textwidth]{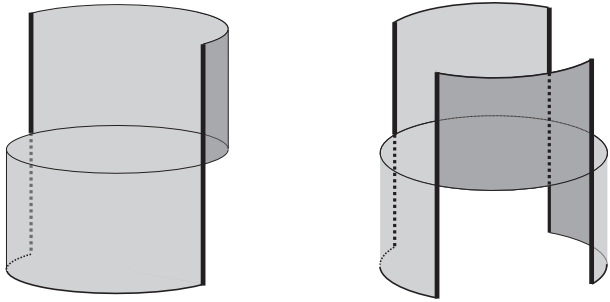}
\caption{$Q'_1(\pi/2)$ and $Q'_2(\pi/2)$.\label{q1'&q2'}}
\end{figure}

\begin{figure}
\centering
\includegraphics[width=0.6\textwidth]{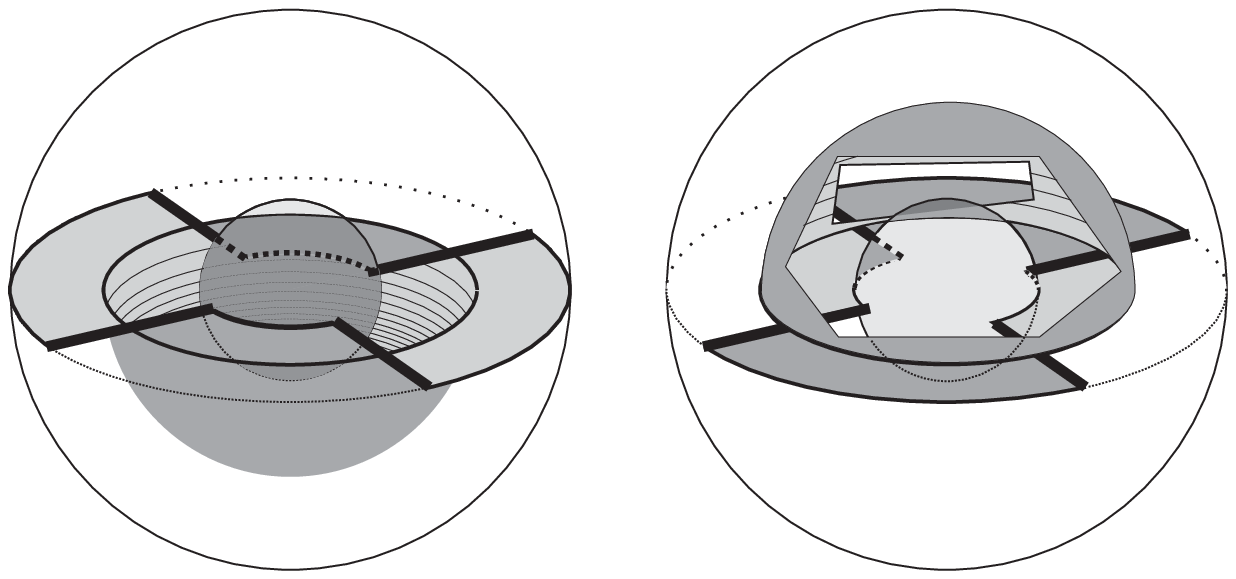}
\caption{$Q'_2(\pi/2)$ and $Q'_2(3\pi/2)$ 
(the latter with viewports),
as they appear in in 
$\Nb{{S^3}}{{S^2}}\protect\iso \protect\Cext\times\clint{0}{\pi}$.
\label{spherical q2'}}
\end{figure}

By further inspection, for
$0<\theta<2\pi$, $\theta\ne\pi$, 
$Q_n(\theta)$ is isotopic
(by a piecewise-smooth isotopy fixing 
$\Bd Q_n(\theta)\cup (Q_n\cap %
\Cext\times\{\theta-\lfloor \theta/\pi \rfloor\pi\})$
pointwise) to a piecewise-smooth $4n$--gonal
$2$--disk $Q'_n(\theta)$ 
so situated that 
$Q'_n(\theta)\cap 
(\Cext\times %
\clint{0}{\theta-\lfloor \theta/\pi \rfloor\pi})$ 
and 
$Q'_n(\theta)\cap %
(\Cext\times %
[\theta-\lfloor \theta/\pi \rfloor\pi,\pi])$ 
are both piecewise-smooth $4n$--gonal $2$--disks,
while $Q'_n(\theta)\cap %
(\Cext\times \{\theta-\lfloor \theta/\pi \rfloor\pi\})$ 
is a (smooth) $2$--disk in 
$\Cext\times\{\theta-\lfloor \theta/\pi \rfloor\pi\}$ 
naturally endowed with the structure of a $2n$--gon.
(The cases $n=1, 2$, $\theta=\pi/2$, are pictured in 
Fig.~\ref{q1'&q2'}.
Versions of $Q'_2(\pi/2)$ and $Q'_2(3\pi/2)$
in $\Nb{{S^3}}{{S^2}}\iso \Cext\times\clint{0}{\pi}$ 
are pictured in Fig.~\ref{spherical q2'}.)

\label{f-good neighborhoods}
Let $L \sub S^3$ be a link, 
$f\from \diff{S^3}{L}\to S^1$ a Morse map, 
and $\psi$ an \bydef{$n$--star},
in the sense of \cite{constqp5},
on a smooth page $\page{f}{\theta}$:
that is, $\psi$ is the union of $n$ arcs $\a_s$,
pairwise disjoint except for a common endpoint
$*_\psi\in\Int \page{f}{\theta}$, and 
$\Bd \page{f}{\theta}\cap\a_s=\diff{\Bd\a_s}{*_\psi}$
for each $s$).
A regular neighborhood $\Nb{S^3}\psi$ is \mbox{\bydef{$f$-good}}
provided that 
\begin{inparaenum}[(a)]
\item\label{f-good patch}
$\Nb{S^3}\psi \cap \page{f}{\theta}$ is a regular neighborhood
$\Nb{\page{f}{\theta}}\psi$ (and thus an \bydef{$n$--patch} on 
$\page{f}{\theta}$ in the sense of 
\cite{constqp5}:
that is, a $2$--disk naturally endowed with 
the structure of a $2n$--gon whose edges are 
alternately boundary arcs and proper arcs in $\page{f}{\theta}$), 
and 
\item\label{f-good match}
there exists a diffeomorphism 
$h\from (\Cext,\G_{2n}) \to (\Bd\Nb{S^3}\psi, L \cap \Bd\Nb{S^3}\psi)$
with $(f\after h)\restr (\diff{\Cext}{\G_{2n}})=p_n$.
\end{inparaenum}

\begin{lem}\label{good patches exist}
Every neighborhood of $\psi$ in $S^3$ contains an 
$f$-good regular neighborhood.
\end{lem}
\begin{proof} 
First suppose that $L=O$ and $f=o$.
Stereographic projection 
\[
\s\from \diff{S^3}{\{(0,-\imunit)\}} \to\C\times\R 
       \suchthat (z,w) \mapsto(z,\Re(w))/(1+\Im(w))
\]
maps $O$ to $S^1\times\{0\}$ and $\page{o}{\pi/2}$ to
$D^2\times\{0\}$; $\psi_1\isdefinedas \s^{-1}(\clint{0}{1}\times\{0\})$
is a $1$--star in $\page{o}{\pi/2}$, and 
the preimage $\s^{-1}(B)$ of an appropriate ellipsoidal 
$3$--disk $B\sub\C\times\R$ (say, with one focus at $(0,0)$,
center at $(1-\e,0)$ for sufficiently small $\e>0$,
and minor axes much shorter than the major axis),
is an $o$\,-good regular neighborhood $\Nb{S^3}{\psi_1}$.
The smooth map 
\[
(\C\times\R)\cup\{\infty\}\to(\C\times\R)\cup\{\infty\} %
\suchthat (\zeta,t)\mapsto(\zeta^n,t),
\thinspace
\infty\mapsto\infty,
\]
can be modified in a neighborhood of $\infty$ 
so as to induce via $\s^{-1}$ a cyclic branched covering 
$c_n\from S^3\to S^3$ of degree $n$,
branched along $A\isdefinedas \{0\}\times S^1\sub S^3$, with 
$c_n^{-1}(\page{o}{\pi/2})=\page{o}{\pi/2}$; then
$\psi_n\isdefinedas c_n^{-1}(\psi_1)$ is an $n$--star
in $\page{o}{\pi/2}$, and (if minimal care has been taken) 
$c_n^{-1}(\Nb{S^3}{\psi_1})$
is an $o$\,-good regular neighborhood $\Nb{S^3}{\psi_n}$.
(The cases $n=1,2,3$ are pictured in 
Fig.~\ref{o-good neighborhoods}.) 
Clearly any neighborhood of $\psi_n$ in $S^3$
contains an $o$\,-good regular neighborhood of the type
just constructed.

\begin{figure}
\centering
\includegraphics[width=0.7\textwidth]{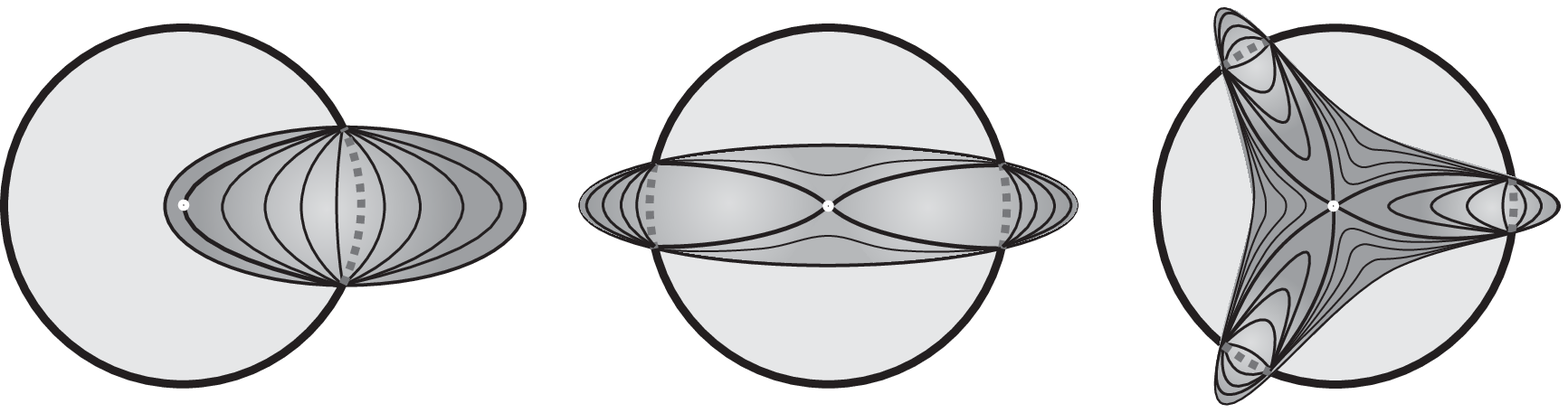}
\caption{Some level sets of 
$o\protect\restr (\diff{\Bd\Nb{S^3}{\psi_n}}{O})$
for $n=1,2,3$.
\label{o-good neighborhoods}}
\end{figure}

The general case follows immediately, upon observing that,
for any link $L$, 
Morse map $f\from \diff{S^3}{L}\to S^1$, 
and spanning surface $\page{f}{\theta}$,
if $\psi\sub \page{f}{\theta}$ is an $n$--star on $\page{f}{\theta}$, 
then there is a diffeomorphism $h\from M\to h(M)$ 
from a (not necessarily $f$-good) neighborhood $M$ of $\psi$ in $S^3$ 
to a neighborhood $h(M)$ of $\psi_n$ in $S^3$
such that $h(M\cap \page{f}{\theta})=h(M)\cap \page{o}{\pi/2}$,
and a diffeomorphism $k\from (S^1,\exp(\imunit\theta))\to (S^1,\imunit)$
such that $k\after f\restr (\diff{M}{L})=o\after h\restr (\diff{M}{L})$.
(A ``side view'' of part of one such $M$ is pictured in 
Fig.~\ref{sideview}.)  
\begin{figure}
\centering
\includegraphics[width=0.4\textwidth]{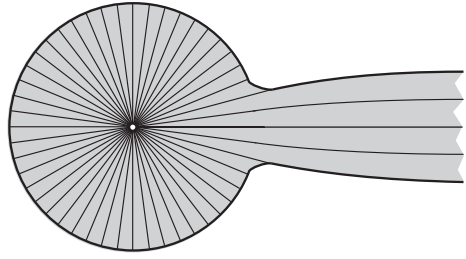}
\caption{A ``side view'' (indicating some level sets
of $f$) of the part of the neighborhood 
$M$ near an endpoint of $\psi$.\label{sideview}}
\end{figure}
\end{proof}

\begin{construction}[Murasugi sum of Morse maps]
\label{M. sum of Morse maps construction}
\hypertarget{M. sum of Morse maps hyperconstruction}{}%
Let $L_s \sub S^3_s$ be a link, 
$f_s\from \diff{S^3_s}{L_s}\to S^1$ a Morse map,
$\exp(\imunit\Theta)\in {S^1}\setminus{\left(f_0(\crit{f_0}) %
\cup f_1(\crit{f_1}\right)\cup\{0,\pi\})}$, 
$\psi_{n,s}\sub \page{f_s}{\Theta}$ an $n$--star
for $s=0,1$ and $n>0$. 
Let $N_s\isdefinedas \Nb{S^3_s}{\psi_{n,s}}$ be $f_s$-good,
$E_s\isdefinedas \Ext{S^3_s}{\psi_{n,s}}$,
$d_s\from E_s\to D^3$ a diffeomorphism.
Let 
$h_s\from (\Cext,\G_{2n}) \to (\Bd N_s, L \cap \Bd N_s)$
be a diffeomorphism such that 
$(f_s\after h_s)\restr (\diff{\Cext}{\G_{2n}})=p_n$.
For a fixed $\zeta\in \G_n$, define $\equiv_\zeta$
by $(z,s)\equiv_\zeta h_s(\zeta^s z,s)$ 
($z\in \Cext, s=0,1$).  
The identification space 
\[
\Sigma(N_0,N_1,\zeta)\isdefinedas  %
\big(E_0 \sqcup %
\thinspace %
(\Cext\times\clint{0}{\pi}) %
\thinspace %
\sqcup E_1 )\big){/}{\equiv}_\zeta
\]
has a natural piecewise-smooth structure,
imposed on it by the identification map $\Pi$,
with respect to which the $1$--submanifold and map
\begin{gather*}
L(N_0,N_1,\zeta)\isdefinedas (L_0\cap E_0 \sqcup %
\G_{2n}\times\clint{0}{\pi}\sqcup L_1\cap E_1){/}%
         {\equiv}_\zeta%
,\\
f(N_0,N_1,\zeta)\isdefinedas %
((f_0 \sqcup P_n \sqcup f_1){/}{\equiv}_\zeta)\restr 
(\diff{\Sigma(N_0,N_1,\zeta)}{L(N_0,N_1,\zeta)}) 
\end{gather*}
are smooth where this is meaningful (i.e., in the complement
of the identification locus $\Pi(\Cext\times\{0,\pi\})$,
along which $\Sigma(N_0,N_1,\zeta)$ is itself 
\foreign{a priori} only piecewise-smooth).
It is not difficult to give $\Sigma(N_0,N_1,\zeta)$ 
a smooth structure everywhere, in which 
$L(N_0,N_1,\zeta)$ and $f(N_0,N_1,\zeta)$ 
are everywhere smooth.
The smoothing can be done very naturally off
$\Pi((\G_{2n}\cup\{0,\infty\})\times\{0,\pi\})$, 
using the fact that by construction 
(see Lemmas~\ref{description of pn and Pn} and
\ref{description of Qn}) 
the function 
$(\|d_0\|^2-1) \sqcup \pr_2 \sqcup (\pi+1-\|d_1\|^2){/}{\equiv}_\zeta$
and 
the map $f(N_0,N_1,\zeta)$ 
are piecewise-smoothly transverse (in an evident sense) on 
$\Pi((\diff{\Cext}{\G_{2n}\cup\{0,\infty\})})\times\{0,\pi\})$.
A somewhat less natural, but not difficult,
construction smooths $\Sigma(N_0,N_1,\zeta)$ on 
$\Pi(\{0,\infty\}\times\{0,\pi\})$ in such a way 
as to make $f(N_0,N_1,\zeta)$ smooth there also
(in the process, 
$(\|d_0\|^2-1) \sqcup \pr_2 \sqcup (\pi+1-\|d_1\|^2){/}{\equiv}_\zeta$
may be forced to be not smooth at those points).
Nor is there is any difficulty in 
smoothing $\Sigma(N_0,N_1,\zeta)$ 
on $\Pi(\G_{2n}\times\{0,\pi\})$.  
Further details will be suppressed.

When $\Sigma(N_0,N_1,\zeta)$, with the smooth structure just
constructed, is identified with $S^3$, 
$L(N_0,N_1,\zeta)$ (resp., $f(N_0,N_1,\zeta)$) 
will be called a \bydef{$2n$--gonal Murasugi sum} 
of $L_0$ and $L_1$ (resp., of $f_0$ and $f_1$) 
and denoted by $L_0\plumb{(N_0,N_1,\zeta)}L_1$ 
(resp., $f_0\plumb{(N_0,N_1,\zeta)}f_1$) or simply 
by $L_0\plumb{}L_1$ (resp., $f_0\plumb{}f_1$).
(A $2$--gonal Murasugi sum of links is simply a 
connected sum, and a $2$--gonal Murasugi sum of Morse
maps $f_0$ and $f_1$ may be denoted $f_0\connsum f_1$.
A more detailed description of the construction
in the case of connected sums, which goes more smoothly 
than the general case, is given in \cite{P-R-W}.)
\end{construction}

\begin{unproved}{thm}\label{properties of Murasugi sum of maps}
Let $L_s \sub S^3_s$ be a link and 
$f_s\from \diff{S^3_s}{L_s}\to S^1$ a Morse map for $s=0,1$.
\begin{inparaenum}
\item\label{sum of links is a link}
Any Murasugi sum $L_0\plumb{(N_0,N_1,\zeta)}L_1$ is a link.
\item\label{sum of maps is a map}
Any Murasugi sum $f_0\plumb{(N_0,N_1,\zeta)}f_1$ is a Morse map,
and its critical points 
$\crit{f_0\plumb{(N_0,N_1,\zeta)}f_1})=\crit{f_0}\cup\crit{f_1}$
have indices that are inherited unchanged
from those of $f_0$ and $f_1$; in particular,
if $f_0$ and $f_1$ are moderate, 
then $f_0\plumb{(N_0,N_1,\zeta)}f_1$ is moderate.
\end{inparaenum}
\end{unproved}

With $\Theta$ as in Construction~\ref{M. sum of Morse maps construction}, 
the Seifert surface 
$\page{f_0\plumb{(N_0,N_1,\zeta)}f_1}{\Theta}$ 
is isotopic to 
$\Pi((\page{f_0}{\Theta}\cap E_0)\sqcup Q_n(\Theta) %
 \sqcup (\page{f_1}{\Theta}\cap E_1))$ (up to smoothing), and thus
piecewise-smoothly isotopic to
\begin{multline*}
\Pi(\page{f_0}{\Theta} \cap E_0 %
   \sqcup Q'_n(\Theta) %
   \sqcup \page{f_1}{\Theta}\cap E_1) = \\
                  \Pi(\page{f_0}{\Theta}\cap E_0\sqcup %
                  Q'_n(\Theta)\cap %
                  \Cext\times 
                  \clint{0}{\Theta-\lfloor \Theta/\pi \rfloor\pi}) 
                  \medspace\cup\qquad\quad\\
\Pi(Q'_n(\Theta)\cap %
                  \Cext\times 
                  \clint{\Theta-\lfloor \Theta/\pi \rfloor\pi}{\pi}
                  \sqcup \page{f_1}{\Theta}\cap E_1)\isdefinedas S'_0\cup S'_1
\end{multline*}
where $S'_s$ is piecewise-smoothly isotopic to $\page{f_s}{\Theta}$
by an isotopy carrying the $2$--disk $S'_0\cap S'_1$ (with
its natural structure of $2n$--gon, noted after 
Lemma~\ref{description of Qn}) to the $n$--patch
$N_s\cap \page{f_s}{\Theta}$.
Thus, for a suitable 
diffeomorphism
$H\from N_0\cap \page{f_0}{\Theta}\to N_1\cap \page{f_1}{\Theta}$
(determined by $\zeta$, up to isotopy), 
$\page{f_0\plumb{(N_0,N_1,\zeta)}f_1}{\Theta}$ 
is a $2n$--gonal Murasugi sum 
$\page{f_0}{\Theta}\plumb{H}\page{f_1}{\Theta}$
as described in \cite{constqp5}
(see also the primary sources 
\cite{Murasugi:plumbing,Stallings,Gabai:Murasugi1}).

On the level of Seifert surfaces, 
a $2$--gonal Murasugi sum $S_0\plumb{H}S_1$
is the same as a boundary-connected sum 
$S_0\bdconnsum S_1$.  After boundary-connected sum, 
the most commonly encountered case of Murasugi sum is
$4$--gonal, and the most familiar and probably most
useful $4$--gonal Murasugi sums are \bydef{annulus plumbings},
as described in \cite{espaliers} (see also the primary
sources \cite{Siebenmann:rational,Gabai:arborescent},
as well as 
\cite{Sakuma:specialarborescent}).

Specifically, 
on any \label{annular plumbing} annulus $\AKn{K}{n}$,
let $\g(K)$ be a proper arc joining the two components
of $\Bd\AKn{K}{n}$,
so that $\Nb{\AKn{K}{n}}{\g(K)}\whichdefines C_1$ 
is a $2$--patch (though $\g(K)$, as given, is not a $2$--star).
Let $S$ be a Seifert surface $S$, 
$\a\sub S$ a proper arc,
$C_0\sub S$ a $2$--patch 
with $\a\sub\Bd C_0$ (respecting orientation).
Let $H\from (C_0,\a)\to (C_1,C_1\cap K)$ be a diffeomorphism.
Each of $\g(K)$, $C_0$, and $H$ is unique up to isotopy,
so the $4$--gonal Murasugi sum $S\plumb{H} \AKn{K}{n}$ 
depends (up to isotopy) only on $\a$, 
and there is no abuse in denoting it by 
$S\plumb{\a} \AKn{K}{n}$.  When $S=\AKn{K'}{n'}$ is also an 
annulus, it is slightly abusive---but handy---to write 
simply $\AKn{K'}{n'}\plumb{}\AKn{K}{n}$, with the understanding 
that $\a=\g(K')$.

\begin{example}\label{connected sum}
$\page{f_0\connsum f_1}{\theta}$ is a 
boundary-connected sum $\page{f_0}{\theta}\bdconnsum \page{f_1}{\theta}$,
for any $f_0$ and $f_1$. 
\end{example}

\begin{example}\label{u*u} 
As pictured in Fig.~\ref{uround} (redrawn piecewise-smoothly
in Fig.~\ref{piecewise-smooth U})
the large Seifert surface of $u$ is $\AKn{O}{0}$.
Of course the small Seifert surface of $u$ is two $2$--disks.  
By the construction and Thm.~\ref{properties of Murasugi sum of maps},
there is a $4$--gonal Murasugi sum $u\plumb{}u$ 
with $4$ critical points and $4$ critical values.
Each non-singular page of $u\plumb{}u$ is either a once-punctured 
torus bounded by an unknot or a $2$--disk bounded by an unknot
(see Fig.~\ref{piecewise-smooth U*U}); note that, 
on the level of Murasugi sums of surfaces, a $2$--disk
is not a Murasugi sum of two pairs of two $2$--disks.
The Morse map $u\plumb{}u$ may be modified in the standard way
to yield a self-indexed Morse map $w$, with small surface 
a $2$--disk and large surface a once-punctured surface of 
genus $2$; $w$ is not a Murasugi sum of Morse maps.  
\begin{figure}
\centering
\includegraphics[width=0.75\textwidth]{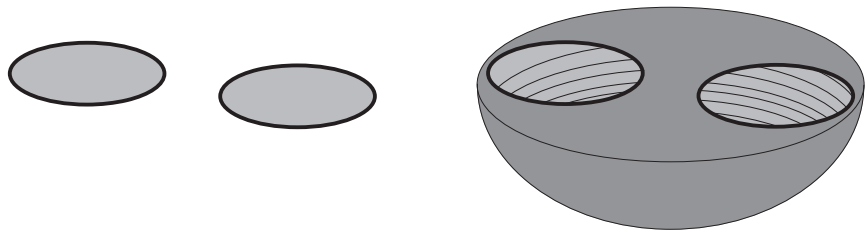}
\caption{Piecewise-smooth Seifert surfaces for $U$ 
(two $2$-disks; an annulus)
isotopic to those in Fig. \protect\ref{uround}.
\label{piecewise-smooth U}}
\end{figure}
\begin{figure}
\centering
\includegraphics[width=0.75\textwidth]{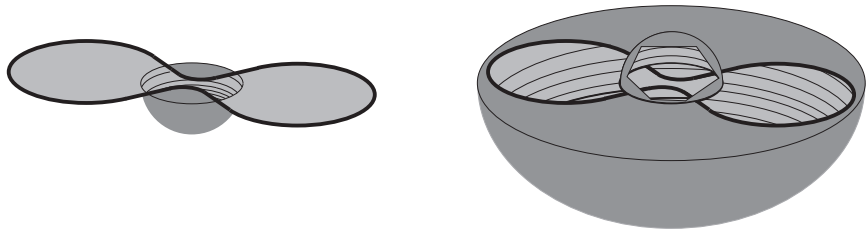}
\caption{Piecewise-smooth Seifert surfaces for 
$U\plumb{}U\protect\iso O$ 
(a $2$-disk; a punctured torus, with viewport).
\label{piecewise-smooth U*U}}
\end{figure}
\end{example}

\begin{example}\label{u*o{2,2}}
Since $o\{2,2\}$ is a fibration, the $4$--gonal Murasugi 
sum $u\plumb{}o\{2,2\}$ is moderate and self-indexed,
with binding isotopic to $O$; in fact, $u\plumb{}o\{2,2\}$ 
is isotopic to $o_1$.  Note that the large and small Seifert 
surfaces of $u\plumb{}o\{2,2\}$ are isotopic to the appropriate
smooth pages of $u\plumb{}u$.
\end{example}

\begin{example}\label{special arborescent links, baskets}
As in \cite{constqp5}, a \bydef{basket} is a Seifert surface 
of the form $D^2\plumb{a_1}\AKn{O}{k_1}\dotsm\plumb{a_n} \AKn{O}{k_n}$,
where $\a_1,\dots,\a_n\sub D^2$ are proper arcs with pairwise
disjoint endpoints.  By Cor.~\ref{MN Bd A{O,n} is 2, except} and
Construction~\ref{M. sum of Morse maps construction},
if $S=D^2\plumb{a_1}\AKn{O}{k_1}\dotsm\plumb{a_n} \AKn{O}{k_n}$ 
is a basket, then $\MN\Bd S\le 2m$, where 
$m\isdefinedas\card{\{j:1\le j\le n \text{ and } |k_n|\ne 2\}}$.
In particular, by \cite{espaliers}, if $L$ is a 
\bydef{special arborescent link} in the sense of 
Sakuma \cite{Sakuma:specialarborescent} (that is, 
if $L$ bounds an arborescent plumbing of annuli $\AKn{O}{k}$),
then $L$ bounds a basket (which is isotopic to such an arborescent
plumbing) and $\MN\Bd S\le 2m$, where again $m$ is the 
number of annular plumbands which are not Hopf annuli.
\end{example}

As Example~\ref{u*o{2,2}} shows, if no restriction is placed on
the Seifert surfaces along which a Murasugi sum $L_0\plumb{}L_1$ 
is formed, then it can easily happen that 
$\MN(L_0)+\MN(L_1) < \MN(L_0\plumb{}L_1)$.
Much worse is true: as observed in \cite{Hirasawa:Murasugi},
any knot $K$ bounds a Seifert surface $S$ such that 
$S=S_0\plumb{}S_1$ is a Murasugi sum 
and $K_0\isdefinedas \Bd S_0$, $K_1\isdefinedas \Bd S_1$ are unknots; 
so $\MN(K_0\plumb{}K_1)-(\MN(K_0)+\MN(K_1))$ can 
be arbitrarily large, since \cite{P-R-W} for every $m$ 
there is a knot $K$ with $\MN(K)>m$.
However, Thm.~\ref{properties of Murasugi sum of maps}
does lead immediately to the following generalization
of the inequality 
\thetag{\ref{subadditivity over connected sum:intro}}
stated in Section~\ref{introduction}.

\begin{unproved}{cor}\label{subadditivity over Murasugi sum--cor}
If $f_s\from \diff{S^3_s}{L_s}\to S^1$ are minimal Morse maps
$(s=0, 1)$, then 
\begin{equation}
\MN(L_0\plumb{(N_0,N_1,\zeta)}L_1)\le \MN(L_0)+\MN(L_1)
\tag{$\ast$}
\label{subadditivity over Murasugi sum:proof}%
\end{equation}%
for any Murasugi sum $f_0\plumb{(N_0,N_1,\zeta)}f_1$.
\end{unproved}

One special, very simple case of Murasugi sum 
is so useful for applications that we draw
explicit attention to it.

\begin{construction}[cutting]
\label{cutting construction}
\hypertarget{cutting hyperconstruction}{}
Let $f\from \diff{S^3}{L}\to S^1$ be a Morse map.
Let $\a\sub\page{f}{\theta}$ be a proper arc on a smooth
page of $f$.  Among the smooth pages of $f\plumb{\a}u$ 
are surfaces isotopic to $\page{f}{\theta}\plumb{\a}\AKn{O}{0}$ 
(that is, $\page{f}{\theta}$ with an untwisted, unknotted
$1$--handle attached ``running above $\a$'') and 
surfaces isotopic to 
$\cut{\page{f}{\theta}}{\a}\isdefinedas\Ext{\page{f}{\theta}}{\a}$
(that is, $\page{f}{\theta}$ ``cut along $\a$'').
Let the Morse map $\cut{f}{\a}\isdefinedas f\plumb{\a}u$ 
be called \bydef{$f$ cut along $\a$}; denote the binding of $\cut{f}{\a}$ 
by $\cut{L}{\a}$ and call it \bydef{$L$ cut along $\a$}.
Of course $\cut{L}{\a}$, $\Bd(\page{f}{\theta}\plumb{\a}\AKn{O}{0})$,
and $\Bd(\cut{\page{f}{\theta}}{\a})$ are mutually isotopic.
\end{construction}

\begin{unproved}{lem}\label{effect of cutting on MN}
If $f\from \diff{S^3}{L}\to S^1$ is a Morse map,
and $\a\sub\page{f}{\theta}$ is a proper arc on any smooth
page of $f$, then 
$\MN(\cut{L}{\a})\le \MN(L)+2$.
\end{unproved}

Given a Seifert surface $S$ and a proper arc $\a\sub S$, 
let $\cross{\a}\sub S$ be a proper arc such that 
$\cross{\a}\sub \Nb{S}{\a}$ and $\cross{\a}$ has 
exactly one point of intersection with $\a$, at which
the arcs are transverse.  Let $\e$ be
$1$ if the arc of $\Bd S\cap \Nb{S}{\a}$ is oriented
from $\a$ to $\cross{\a}$, $-1$ in the contrary case.

\begin{lem}\label{undoing a cut} The links $\Bd S$ and 
$\Bd(S\plumb{\a}\AKn{O}{0}\plumb{\cross{\a}}\AKn{O}{-\e})$
are isotopic.
\end{lem}
\begin{proof} See Fig.~\ref{undo cut}.
\begin{figure}
\centering
\includegraphics[width=0.6\textwidth]{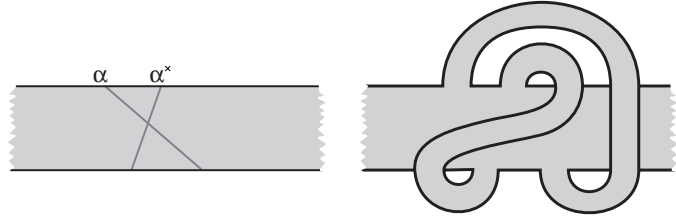}
\caption{The Seifert surfaces $S$ (left)
and $S\plumb{\a}\AKn{O}{0}\plumb{\cross{\a}}\AKn{O}{-1}$ 
(right) have isotopic boundaries.
\label{undo cut}}
\end{figure}
\end{proof}

\begin{hist}\label{Goda's priority}
After the second author had prepared \cite{Rudolph:Msums}, 
the first author brought to his attention work of Goda 
\cite{Goda:suture-sum,Goda:handle-number} which immediately 
implies \thetag{\ref{subadditivity over Murasugi sum:proof}},
although Goda's results are stated in the language of 
``handle number of a Seifert surface'' and 
``Murasugi sum of Seifert surfaces''
rather than that of ``Morse--Novikov number of a link'' 
and ``Murasugi sum of Morse maps'', and Goda's proofs 
correspondingly use the techniques of sutured manifolds 
\cite{Gabai:foliations-genera} and 
$C$-product decompositions \cite{Gabai:detect} 
rather than those of Morse maps.  We have given the 
proof using Morse maps for the sake of variety.
\end{hist}

\section{Morse maps and free Seifert surfaces}\label{freeness}

A Seifert surface $S\sub S^3$ is called 
\bydef{free}\label{definition of free}
if and only if $S$ is connected and the group $\pi_1(\diff{S^3}{S})$ 
is free; alternatively, $S$ is free if and only if 
$\diff{S^3}{\Int\Nb{S^3}{S}}$ is a handlebody 
$\hb{g}$ (necessarily of genus $g=b_1(S)$,
where for any $X$, we write $b_1(X)\isdefinedas\rank H_1(X;\Z)$).
A Morse map $f$ is \bydef{free} provided that $f$ is connected
and every Seifert surface $\page{f}{\theta}$ is free.
It is well known (and obvious) that, if $f$ is a fibration,
then $f$ is free.  
The following proposition is Lemma~4.2 of \cite{P-R-W}.

\begin{unproved}{prop}
If $f$ is a moderate self-indexed Morse map, then 
a large Seifert surface of $f$ is free.
\end{unproved}

\begin{unproved}{cor} 
If $f$ is a moderate self-indexed Morse map
\textup{(}in particular, if $f$ is a self-indexed 
minimal Morse map\textup{)}, then $f$ is free if 
and only if a small Seifert surface of $f$ is free.
\end{unproved}

\begin{thm} If $S$ is a free Seifert surface,
then there is a free moderate self-indexed Morse
map $f\from \diff{S^3}{\Bd S}\to S^1$, with 
exactly $2b_1(S)$ critical points, such that
$S$ is the small Seifert surface of $f$.
\end{thm}
\begin{proof}
Given a connected Seifert surface $S$
with $b_1(S)\whichdefines n$, enlarge it to 
a Seifert surface $S^{+}$ as follows.
Let $\a_1(S),\dots,\a_n(S)\sub S$ be pairwise disjoint proper arcs 
such that $\cut{\cut{\cut{\cut{S}{\a_1(S)}}{\a_2(S)}}{\dotsm}}{\a_n(S)}$
is a $2$--disk.
There is a handlebody $\hb{n}\sub S^3$ such that
$S\sub\Bd\hb{n}$ and $\hb{n}$ is the trace of an isotopy
(rel.~$\Bd S$) from $S$ to $\diff{\Bd\hb{n}}{\Int S}$.
The trace of $\a_s(S)$ by a suitable such isotopy is a meridional 
$2$--disk of $\hb{n}$ with unknotted boundary $\g_s(S)$,
and $\Nb{\Bd\hb{n}}{\g_s(S)}$ is an annulus 
$\AKn{\g_s(S)}{0}$.  Evidently
\[
S^{+}\isdefinedas S\cup\bigcup_s \AKn{\g_s(S)}{0} = 
S\plumb{\a_1(S)}\AKn{\g_1(S)}{0}\dotsm\plumb{\a_n(S)}\AKn{\g_n(S)}{0}
\]
is a subsurface of $\Bd\hb{n}$ bounded by an unknot.
If also $S$ is free, so that $\hb{n}$ is Heegaard,
then by Waldhausen's uniqueness theorem \cite{Waldhausen1968}
up to isotopy the pair $(\Bd\hb{n},S^{+})$ depends only
on $n$ and is otherwise independent of $S$. 
In particular, if $T$ is a fiber surface
with $b_1(T)=n$ (which always exists: e.g., a boundary-connected 
sum of $n$ Hopf annuli), then $T^{+}$ is isotopic to $S^{+}$.
Now, $T$ is a large Seifert surface of a fibration $f_T$;
by Construction~\ref{M. sum of Morse maps construction},
Thm.~\ref{properties of Murasugi sum of maps}, 
and Construction~\ref{cutting construction}, $T^{+}$ is a large 
Seifert surface of a self-indexed moderate Morse map 
$\cut{\cut{\cut{\cut{f_T}{\a_1(T)}}{\a_2(T)}}{\dotsm}}{\a_n(T)}$
with $2n$ critical points.  Let $\b_1,\dots,\b_n\sub T^{+}$
be the arcs (which without loss of generality can be chosen to 
be proper) which are the images of suitably chosen
$\cross{\a_1(S)}, \dots, \cross{\a_s(S)}\sub S\sub S^{+}$ 
by some isotopy carrying $(\Bd\hb{n},S^{+})$ onto 
$(\Bd\hb{n},T^{+})$.  
(The situation is illustrated in 
Fig.~\ref{free genus estimate figure}.)
\begin{figure}
\centering
\includegraphics[width=.74\textwidth]{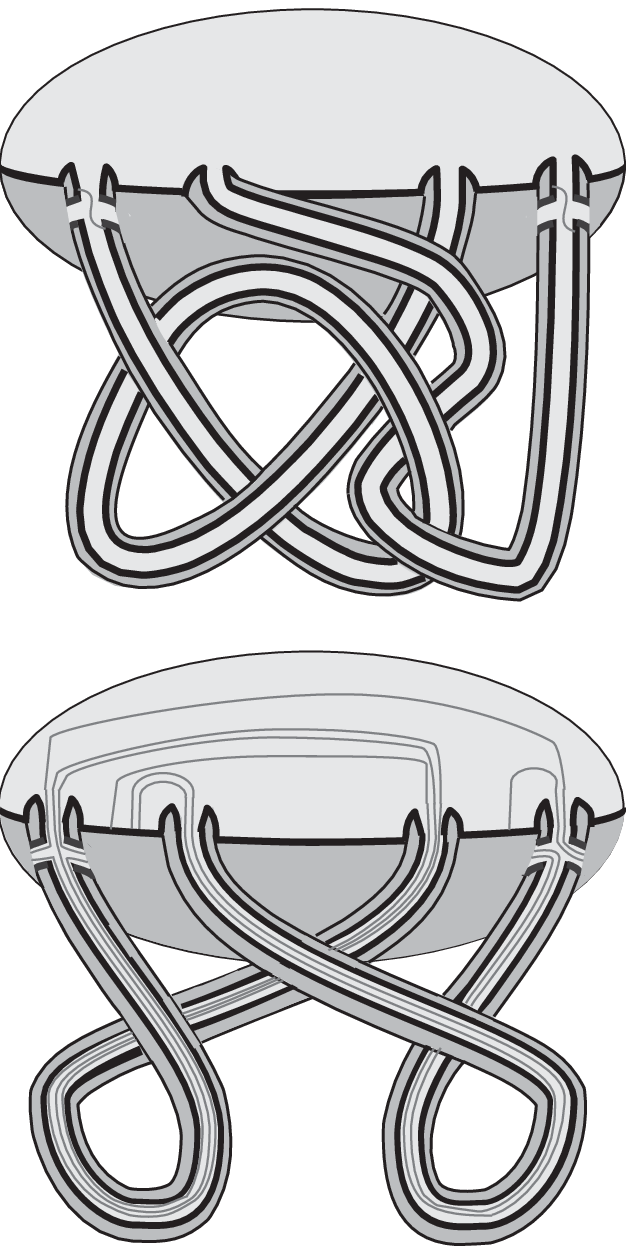}
\caption{A pair of isotopic triples
$(\Bd\hb{2},S^{+},(\cross{\a_1(S)},\cross{\a_2(S)}))$ 
and $(\Bd\hb{2},T^{+},(\b_1,\b_2))$ constructed
from a free Seifert surface $S$ and a fiber surface $T$.
\label{free genus estimate figure}}
\end{figure}
Again by Construction~\ref{M. sum of Morse maps construction} 
and Thm.~\ref{properties of Murasugi sum of maps}, 
\[
(\cut{\cut{\cut{\cut{f_T}{\a_1(T)}}{\a_2(T)}}{\dotsm}}{\a_n(T)})%
\plumb{\b_1}o\{2,-2\e_1\}\dotsm\plumb{\b_n}o\{2,-2\e_n\}
\]
is a self-indexed moderate Morse map with $2n$ critical points,
and by Lemma~\ref{undoing a cut} its binding 
$\Bd(T^{+}\plumb{\b_1}O\{2,-2\e_1\}\dotsm\plumb{\b_n}O\{2,-2\e_n\})$
is isotopic to $\Bd S$, while evidently its small Seifert
surface is isotopic to $S$.
\end{proof}

The \bydef{free genus} of a knot $K$ is 
$\freegenus{K}\isdefinedas \min\{g(S)\Suchthat 
K=\Bd S$, $S$ is a free Seifert surface$\}$.
More generally, for any link $L$ let 
$\freerank{L}\isdefinedas \min\{n\Suchthat
L=\Bd S$, $S$ is a free Seifert surface$\}$; 
so, for instance, $\freerank{K}=2\freegenus{K}$ if $L=K$ is a knot, 
and $\freerank{L}=1$ if and only if $L=\Bd \AKn{O}{n}$ 
for some $n\in\Z$, $n\ne 0$.

\begin{unproved}{free genus estimate}
\label{MN bounded by four times the free genus} 
\hypertarget{free genus estimate theorem}%
{For any link $L$, $\MN(L)\le 2\freerank{L}$.  
In particular, for any knot $K$, $\MN(K)\le 4\freegenus{K}$.}
\end{unproved}

\begin{remarks} 
\begin{inparaenum}
\item
Let $g(K)$, as usual, denote the genus of a knot $K$.
Of course $\freegenus{K}\ge g(K)$ for every knot $K$.
For many knots $K$ (satisfying a suitable condition on
the Alexander polynomial $\D_K(t)$), it follows 
from \cite{P-R-W} that $\MN(K)\ge 2g(K)$.  We know no 
example of a knot $K$ for which it can be shown that 
$\MN(K)>2g(K)$.  There are knots $K$ with $g(K)=1$ 
and $\freegenus{K}$ arbitrarily large \cite{Moriah:free,Livingston:free}.
\item
Clearly the 
{\hyperlink{free genus estimate}{Free Genus Estimate}}
can be formally strengthened to aver that 
$\MNfree(L)\le 2\freerank{L}$ for any link $L$
and $\MNfree(K)\le 4\freegenus{K}$ for any knot $K$, 
where the \bydef{free Morse--Novikov number} $\MNfree(L)$ 
of $L$ is the minimum possible number of critical points of a free 
Morse map with binding $L$.  There exist knots and links with minimal 
Morse maps which are not free.  It would be interesting to 
know whether there exists $L$ with $\MN(L)<\MNfree(L)$.
\end{inparaenum}
\end{remarks}

\begin{hist}\label{history of free genus}
Neuwirth \cite{Neuwirth:algebraic} was 
apparently the first to consider, in effect, the
notion of ``free Seifert surface'': in a footnote, he 
called a Seifert surface $S$ ``algebraically knotted'' 
if $\pi_1( \diff{S^3}{S})$ is not free.  
Neuwirth's language was adapted by 
Murasugi \cite{Murasugi:commutator}, 
Lyon \cite{Lyon:1972}, and others.
By 1976, Problem~1.20 (attributed to Giffen and Siebenmann)
in Kirby's problem list \cite{Kirby:problems}
uses the phrase ``free Seifert surface'', 
points out that $S$ is free if and only if $\diff{S^3}{S}$ 
is an open handlebody,
apparently introduces the terminology ``free genus'',
and includes the imperative ``Relate the free genus to 
other invariants of knots.''
A number of later authors
(e.g., Moriah \cite{Moriah:free},
Livingston \cite{Livingston:free},
M.~Kobayashi \& T.~Kobayashi \cite{Kobayashi-Kobayashi})
have studied free Seifert surfaces.
\end{hist}

\section{Morse maps of doubled knots}\label{annuli}

In this section, we estimate the Morse--Novikov number 
of a Whitehead double $\double{K}{m}{\pm}\isdefinedas
\Bd\AKn{K}{k}\plumb{} \AKn{O}{\mp 1}$ of a knot $K$
in terms of various invariants of $K$.

\begin{braid index estimate}
\hypertarget{braid index estimate theorem}%
{For any knot $K$ and all $m\in\Z$, 
$\MN(\double{K}{m}{\pm})\leq 4\bi(K)-2$.}
\end{braid index estimate}
\begin{proof}
Let $n\isdefinedas\bi(K)$.
We may assume that $K$ is a closed $n$-string
$o$-braid.  By Cor.~\ref{splicing doesn't make MN bigger},
$\splicedmap{o}{K}{\framing{k}}$ is a fibration for any $k\in\Z$.
As illustrated in Fig.~\ref{cutting a page}, %
\begin{figure}
\centering
\includegraphics[width=0.7\textwidth]{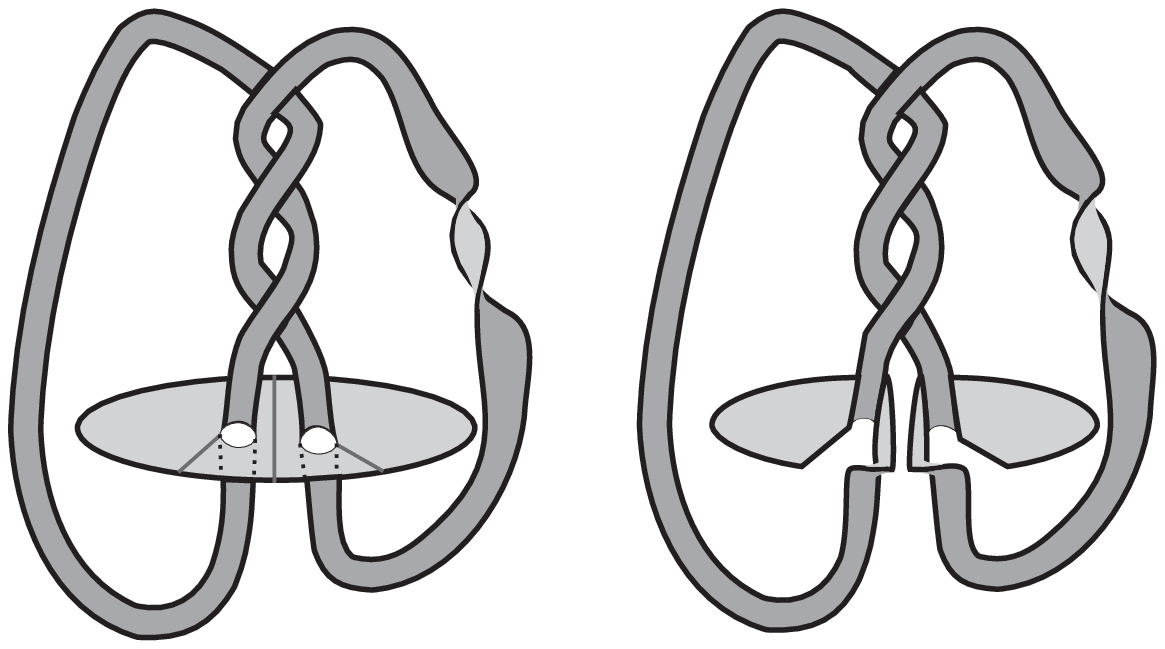}
\caption{\label{cutting a page}}
\end{figure}%
it is easy to take a standard fiber surface $S$ for 
$\splicedmap{o}{K}{\framing{k}}$ (as in 
Construction~\ref{splicing construction})
and find $2n-1$ proper arcs $\a_i\sub S$
such that 
\[
\cut{\cut{\cut{\cut{S}{\a_1}}{\a_2}}{\dotsm}}{\a_{2n-1}}\whichdefines
\AKn{K}{m}
\]
is an annulus; by adjusting $k$, any value of $m$ can be achieved.  
By Lemma~\ref{effect of cutting on MN},
$\MN(\AKn{K}{m})\le 4\bi(K)-2$.  Thus 
$\MN(\double{K}{m}{\pm})=\MN(\AKn{K}{m}\plumb{}\AKn{O}{\mp1})\le 4\bi(K)-2$.
\end{proof}

\begin{wrapping genera estimate}%
\hypertarget{wrapping genera estimate theorem}%
{For any knot $K$ and all $m\in\Z$, 
$\MN(\double{K}{m}{\pm}\leq 2(\wlapgenus{K}+1)\leq 2(\wrapgenus{K}+1)$.}
\end{wrapping genera estimate}%
\begin{proof}
Let $n\isdefinedas\wlapgenus{K}$.
We may assume that $K$ is a closed $1$-string
$o_n$-braid.  By Cor.~\ref{splicing doesn't make MN bigger},
$\MN(\splicedlink{O}{K}{\framing{m-1}}\le 2n$. 
Entirely similarly to the situation pictured in Fig.~\ref{cutting a page}, 
it is easy to take the standard Seifert surface $S$ for 
$\splicedmap{o}{K}{\framing{m-1}}$, and find on it 
a proper arc $\a$ such that $\cut{S}{\a}$ is an annulus
$\AKn{K}{m}$.  By Lemma~\ref{effect of cutting on MN}, 
$\MN(\AKn{K}{m})\le 2n+2$.  Thus 
$\MN(\double{K}{m}{\pm})=\MN(\AKn{K}{m}\plumb{}\AKn{O}{\mp1})\le 
2(\wlapgenus{K}+1)$.  As already noted in 
Cor.~\ref{wlap at most wrap}, $\wlapgenus{K}\le\wrapgenus{K}$.
\end{proof}

\begin{remark} The first author has observed that, for every knot $K$, 
there exists (at least one) $k$ such that 
$\MN(\double{K}{m}{\pm}\leq 2\wrapgenus{K}$
for $m=k-1$ and $m=k+1$.  The proof will be deferred
to another paper; it uses the technique of sutured 
manifolds.
\end{remark}

\begin{crossing number estimate}
\hypertarget{crossing number estimate theorem}%
{For any knot $K$ and all $m\in\Z$, 
$\MN(\double{K}{m}{\pm}\leq 2(c(K)+2)$.}
\end{crossing number estimate}

\begin{proof}
This follows immediately from the estimate 
$\MN(\double{K}{m}{\pm}\leq 2(\wrapgenus{K}+1)$, 
given the trivial observation that $\wrapgenus{K}\le c(K)+1$.
\end{proof}


\begin{thebibliography}{10}

\bibitem{Aristotle}
Aristotle, \emph{Physics}, Book $\Gamma$, Part 5, p.~204a8.

\bibitem{Bodin}
Arnaud Bodin, \emph{Milnor fibration and fibred links at infinity}, 
  Internat. Math. Res. Notices \textbf{11} (1999), 615--621.  \MR{2000b:57046}

\bibitem{Boileau-Orevkov}
Michel Boileau and Stepan~Yu. Orevkov, \emph{Quasipositivit\'e d'une courbe
  analytique dans une boule pseudo-convexe}, C. R. Acad. Sci. Paris
  \textbf{332} (2001), 825--830.  \MR{2002d:32039}

\bibitem{Brown-Crowell:augmentation}
E. M. Brown and R. H. Crowell, 
  \emph{The augmentation subgroup of a link}, 
  J. Math. Mech. \textbf{15} (1966), 1065--1074.
  \href{http://www.ams.org/mathscinet-getitem?mr=33\%20\%234920}%
  {MR 33 \#4920}

\bibitem{Eisenbud-Neumann}
David Eisenbud and Walter Neumann, \emph{Three-dimensional link theory and
  invariants of plane curve singularities}, Princeton University Press,
  Princeton, N.J., 1985, Annals of Mathematics Studies, 110. \MR{87g:57007}

\bibitem{Gabai:Murasugi1}
David Gabai, \emph{The {M}urasugi sum is a natural geometric operation},
  Low-dimensional topology (San Francisco, Calif., 1981), Amer. Math. Soc.,
  Providence, R.I., 1983, pp.~131--143. \MR{85d:57003}

\bibitem{Gabai:foliations-genera}
\bysame, \emph{Foliations and genera of links}, Topology \textbf{23} (1984),
  no.~4, 381--394. \MR{86h:57006}

\bibitem{Gabai:detect}
\bysame, \emph{Detecting fibred links in ${S}\sp 3$}, Comment. Math. Helv.
  \textbf{61} (1986), no.~4, 519--555. \MR{88c:57009}

\bibitem{Gabai:arborescent}
\bysame, \emph{Genera of the arborescent links}, Mem. Amer. Math. Soc.
  \textbf{59} (1986), no.~339, i--viii and 1--98. \MR{87h:57010}

\bibitem{Giroux:fibered}
Emmanuel Giroux, \emph{Fibered links and contact structures}, lectures at
  Research-in-Pairs Workshop on $3$-manifolds and Singularities, Oberwolfach,
  November 19--25 2000.

\bibitem{Giroux:ICM}
\bysame, \emph{{G\'eom\'etrie de contact: de la dimension trois vers les
  dimensions sup\'erieures}}, Proceedings of the {ICM}, Beijing 2002, Vol. $3$,
  2003, \url{http://arxiv.org/abs/math.GT/0305129}, pp.~405--414.

\bibitem{Goda:suture-sum}
Hiroshi Goda, \emph{Heegaard splitting for sutured manifolds and {M}urasugi
  sum}, Osaka J. Math. \textbf{29} (1992), no.~1, 21--40. \MR{93d:57014}

\bibitem{Goda:handle-number}
\bysame, \emph{On handle number of {S}eifert surfaces in ${S}\sp 3$}, Osaka J.
  Math. \textbf{30} (1993), no.~1, 63--80. \MR{93k:57010}

\bibitem{Harer}
John Harer, \emph{How to construct all fibered knots and links}, 
  Topology \textbf{21} (1982), no.~3, 263--280.
  \MR{83e:57007}

\bibitem{Hirasawa:Murasugi}
Mikami Hirasawa, \emph{Counter-intuitive aspects of {M}urasugi sums}, talk at
  {Knots, Links \& $3$-Manifolds---$4$th International Siegen Topology
  Symposium, January 2001}.

\bibitem{Hirasawa-Rudolph:enhanced}
\bysame and Lee Rudolph, 
  \emph{Isolated critical points of maps $\R^4\to\R^2$,
        contact forms on $S^3$, and the enhanced Milnor number}, 
  in preparation (2003).

\bibitem{Katz}
Nicholas Katz, personal communication (June, 2003).

\bibitem{Kirby:problems}
Rob Kirby, \emph{Problems in low dimensional manifold theory}, Algebraic and
  geometric topology (Proc. Sympos. Pure Math., Stanford Univ., Stanford,
  Calif., 1976), Part 2, Amer. Math. Soc., Providence, R.I., 1978, updated
  version at \url{http://math.berkeley.edu/~kirby/problems.ps.gz},
  pp.~273--312. \MR{80g:57002}

\bibitem{Kobayashi-Kobayashi}
Masako Kobayashi and Tsuyoshi Kobayashi, \emph{On canonical genus and free
  genus of knot}, J. Knot Theory Ramifications \textbf{5} (1996), no.~1,
  77--85. \MR{97d:57008}

\bibitem{Lazarev}
A.~Yu. Lazarev, \emph{The {N}ovikov homology in knot theory}, Mat. Zametki
  \textbf{51} (1992), no.~3, 53--57, 144. \MR{93g:57012}

\bibitem{Livingston:free}
Charles Livingston, \emph{The free genus of doubled knots}, Proc. Amer. Math.
  Soc. \textbf{104} (1988), no.~1, 329--333. \MR{89j:57004}

\bibitem{Lyon:1972}
Herbert~C. Lyon, \emph{Knots without unknotted incompressible spanning
  surfaces}, Proc. Amer. Math. Soc. \textbf{35} (1972), 617--620. 
  \href{http://www.ams.org/mathscinet-getitem?mr=46\%20\%232663}{MR 46 \#2663}

\bibitem{Milnor:Morse}
John Milnor, \emph{Morse theory}, Princeton University Press, Princeton, N.J.,
  1963, Annals of Mathematics Studies, No. 51. 
  \href{http://www.ams.org/mathscinet-getitem?mr=29\%20\%23634}{MR 29 \#634}

\bibitem{Milnor:singular-points}
\bysame, \emph{Singular points of complex hypersurfaces}, Princeton University
  Press, Princeton, N.J., 1968, Annals of Mathematics Studies, No. 61. 
  \href{http://www.ams.org/mathscinet-getitem?mr=39\%20\%23969}{MR 39 \#969}

\bibitem{Moriah:free}
Yoav Moriah, \emph{On the free genus of knots}, Proc. Amer. Math. Soc.
  \textbf{99} (1987), no.~2, 373--379. \MR{88a:57012}

\bibitem{Murasugi:plumbing}
Kunio Murasugi, \emph{On a certain subgroup of the group of an alternating
  link}, Amer. J. Math. \textbf{85} (1963), 544--550. 
  \href{http://www.ams.org/mathscinet-getitem?mr=28\%20\%23609}{MR 28 \#609}

\bibitem{Murasugi:commutator}
\bysame, \emph{The commutator subgroups of the alternating knot groups}, Proc.
  Amer. Math. Soc. \textbf{28} (1971), 237--241. 
  \href{http://www.ams.org/mathscinet-getitem?mr=43\%20\%231171}{MR 43 \#1171}

\bibitem{Neumann}
Walter Neumann, 
  \href{http://134.76.163.65/servlet/digbib?template=view.html%
        &id=178422&startpage=451&endpage=496&image-path=%
        http://134.76.176.141/cgi-bin/letgifsfly.cgi%
        &pagenumber=451&imageset-id=4664}
  {\emph{Complex algebraic plane curves via their links at infinity}}, 
    Invent. Math. \textbf{98} (1989), no.~3, 445--489. \MR{91c:57014}

\bibitem{Neumann-Rudolph}
\bysame and Lee Rudolph, 
  \href{http://134.76.163.65/servlet/digbib?template=view.html%
        &id=162055&startpage=423&endpage=454&image-path=%
        http://134.76.176.141/cgi-bin/letgifsfly.cgi%
        &pagenumber=423&imageset-id=4139}
  {\emph{Unfoldings in knot theory}}, Math. Ann.
  \textbf{278} (1987), no.~1-4, 409--439. \MR{89j:57017a};
  \href{http://134.76.163.65/servlet/digbib?template=view.html%
        &id=162251&startpage=355&endpage=358&image-path=%
        http://134.76.176.141/cgi-bin/letgifsfly.cgi%
        &pagenumber=355&imageset-id=4142}{corrigendum},
  Math. Ann. \textbf{282} (1988), no.~2, 349--351. \MR{89j:57017b}

\bibitem{Neuwirth:algebraic}
L.~Neuwirth, \emph{The algebraic determination of the topological type of the
  complement of a knot}, Proc. Amer. Math. Soc. \textbf{12} (1961), 904--906.
\href{http://www.ams.org/mathscinet-getitem?mr=24\%20\%23A3642}{MR 24 \#A3642}

\bibitem{Novikov:multivalued}
S.~P. Novikov, \emph{Multivalued functions and functionals. {A}n analogue of
  the {M}orse theory}, Dokl. Akad. Nauk SSSR \textbf{260} (1981), no.~1,
  31--35. \MR{83a:58025}

\bibitem{P-R-W}
A.~Pajitnov, L.~Rudolph, and C.~Weber, \emph{The {M}orse-{N}ovikov number for
  knots and links} (Russian), Algebra i Analiz \textbf{13} (2001), no. 3,
  105--118; translation in St. Petersburg Math. J. \textbf{13} (2002), no.~3,
  417--426. \MR{2002k:57040}

\bibitem{Pichon}
Anne Pichon, \emph{Germes analytiques r\'eels et d\'ecompositions
  en livre ouvert de la sph\`ere de dimension trois}, in preparation (2003).

\bibitem{embeddings}
Lee Rudolph, \emph{Embeddings of the line in the plane},
  J. Reine Angew. Math. \textbf{337} (1982), 113--118. \MR{84h:14019} 

\bibitem{braided-surfaces}
\bysame, \emph{Braided surfaces and {S}eifert ribbons for closed braids},
  Comment. Math. Helv. \textbf{58} (1983), no.~1, 1--37. \MR{84j:57006}

\bibitem{constqp1}
\bysame, \emph{Constructions of quasipositive knots and links. \textup{I}},
  Knots, braids and singularities (Plans-sur-Bex, 1982), Univ. Gen\`eve,
  Geneva, 1983, pp.~233--245. \MR{86k:57004}

\bibitem{algfuns}
\bysame, \emph{Algebraic functions and closed braids}, 
  Topology \textbf{22} (1983), no.~2, 191--202. \MR{84e:57012}

\bibitem{icp1}
\bysame, \emph{Isolated critical points of mappings from $\mathbb{R}\sp4$ to
  $\mathbb{R}\sp2$ and a natural splitting of the {M}ilnor number of a
  classical fibered link. {{\emph{{I}}}}. {B}asic theory; examples}, Comment.
  Math. Helv. \textbf{62} (1987), no.~4, 630--645,
  \url{http://arxiv.org/abs/math.GT/0203032}. \MR{88k:57009}

\bibitem{constqp3}
\bysame, \emph{Constructions of quasipositive knots and links. 
  \textup{III.} {A}   characterization of quasipositive 
  {S}eifert surfaces}, Topology \textbf{31}
  (1992), no.~2, 231--237. \MR{93g:57014}

\bibitem{totaltan}
\bysame, \emph{Totally tangential links of intersection of complex plane curves
  with round spheres}, Topology '90 (Columbus, OH, 1990), de Gruyter, Berlin,
  1992, pp.~343--349. \MR{94d:57027}

\bibitem{constqp5}
\bysame, \emph{Quasipositive plumbing \textup{(}{C}onstructions of
  quasipositive knots and links. \textup{V)}}, Proc. Amer. Math. Soc.
  \textbf{126} (1998), no.~1, 257--267. \MR{98h:57024}

\bibitem{espaliers}
\bysame, \emph{Hopf plumbing, arborescent {S}eifert surfaces, baskets,
  espaliers, and homogeneous braids}, Topology Appl. \textbf{116} (2001),
  no.~3, 255--277, \url{http://arxiv.org/abs/math.GT/9810095}. \MR{2002g:57018}

\bibitem{Rudolph:Msums}
\bysame, \emph{Murasugi sums of Morse maps to the circle, 
    Morse--Novikov numbers, and free genus of knots}, 
    preprint (2001), \url{http://arxiv.org/abs/math.GT/0108006}.

\bibitem{Sakuma:specialarborescent}
Makoto Sakuma, \emph{Minimal genus {S}eifert surfaces for special arborescent
  links}, Osaka J. Math. \textbf{31} (1994), no.~4, 861--905. \MR{96b:57011}

\bibitem{Siebenmann:rational}
Laurent Siebenmann, \emph{Exercices sur les n{\oe}uds rationnels}, mimeographed
  notes, Orsay, 1975.

\bibitem{Siebenmann:splice}
\bysame, \emph{On vanishing of the Rohlin invariant 
and nonfinitely amphicheiral homology $3$-spheres}, 
Topology Symposium, Siegen 1979 (Proc. Sympos., 
Univ. Siegen, Siegen, 1979), 
Springer, Berlin, 1980, pp. 172--222.  
Lecture Notes in Math., No. 788.  \MR{81k:57011}

\bibitem{Skora}
Richard~K. Skora, 
   \href{http://134.76.163.65/servlet/digbib?template=view.html%
        &id=173890&startpage=179&endpage=194&image-path=%
        http://134.76.176.141/cgi-bin/letgifsfly.cgi%
        &pagenumber=179&imageset-id=4424}
        {\emph{Closed braids in {$3$}-manifolds}}, Math. Zeitschrift
  \textbf{211} (1992), no.~2, 173--187. \MR{93m:57009}

\bibitem{Stallings}
John~R. Stallings, \emph{Constructions of fibred knots and links}, Algebraic
  and geometric topology (Proc. Sympos. Pure Math., Stanford Univ., Stanford,
  Calif., 1976), Part 2, Amer. Math. Soc., Providence, R.I., 1978, pp.~55--60.
  \MR{80e:57004}

\bibitem{Stillwell}
John Stillwell, \emph{The compound crossing number of a knot}, 
 Austral. Math. Soc. Gaz. \textbf{6} (1979), no.~1, 1--10.
  \MR{80e:57005}

\bibitem{Sundheim}
Paul~A. Sundheim, \emph{Reidemeister's theorem for {$3$}-manifolds}, Math.
  Proc. Cambridge Philos. Soc. \textbf{110} (1991), no.~2, 281--292.
  \MR{93b:57012}

\bibitem{Waldhausen1968}
Friedhelm Waldhausen, \emph{Heegaard-{Z}erlegungen der {$3$}-{S}ph\"are},
  Topology \textbf{7} (1968), 195--203. 
  \href{http://www.ams.org/mathscinet-getitem?mr=37\%20\%233576}%
  {MR 37 \#3576}



\bibitem{Wyatt}
William Wyatt, personal communication (June, 2003).

\end{thebibliography}

\providecommand{\bysame}{\leavevmode\hbox to3em{\hrulefill}\medspace}
\makeatletter
\renewcommand\MR[1]{\relax\ifhmode\unskip\spacefactor3000 \space\fi
  \href{http://www.ams.org/mathscinet-getitem?mr=#1}{MR #1}
}
\makeatother

\end{document}